\documentclass[preprint,3p]{elsarticle}

\usepackage{amssymb}
\usepackage{amsthm}

\usepackage{hyperref}
\usepackage{latexsym,amsmath}
\usepackage[shortlabels]{enumitem}

\begin{document}

	\renewcommand{\phi}{\varphi}
	\newcommand{\be}{\begin{equation}}
	\newcommand{\ee}{\end{equation}}
	\newcommand{\ba}{\begin{eqnarray}}
	\newcommand{\ea}{\end{eqnarray}}
	\newcommand{\ban}{\begin{eqnarray*}}
		\newcommand{\ean}{\end{eqnarray*}}
	
	\newcommand{\nul}{{\bf0}}
	\newcommand{\rd}{{\mathbb R}^d}
	\newcommand{\zd}{{\mathbb Z}^{d}}
	\renewcommand{\r}{{\mathbb R}}
	\newcommand{\z} {{\mathbb Z}}
	\newcommand{\cn} {{\mathbb C}}
	\newcommand{\n} {{\mathbb N}}
	
	\newcommand{\ddd}{,\dots,}
	\renewcommand{\lll}{\left(}
	\newcommand{\rrr}{\right)}
	\newcommand{\h}{\widehat}
	\newcommand{\w}{\widetilde}

	\newtheorem{theo}{Theorem}
	\newtheorem{lem}[theo]{Lemma}
	\newtheorem {prop} [theo] {Proposition}
	\newtheorem {coro} [theo] {Corollary}
	\newtheorem {defi} [theo] {Definition}
	\newtheorem {rem} [theo] {Remark}
	\newtheorem {ex} [theo] {Example}

\begin{frontmatter}



\title{From Frame-like Wavelets to Wavelet Frames keeping approximation properties and symmetry}


\author{Krivoshein A.V.\fnref{fn1}}
\ead{krivosheinav@gmail.com, a.krivoshein@spbu.ru}

\fntext[fn1]{The author acknowledges Saint-Petersburg State University for a research grant \# 9.38.198.2015
	and RFBR, research project No. 15-01-05796 a, and Volkswagen Foundation.}

\address{St. Petersburg State University, 7/9 Universitetskaya nab., St. Petersburg, 199034 Russia.}

\begin{abstract}
		For a given symmetric refinable mask obeying the sum rule of order $n$, an explicit method  is suggested for the construction of mutually symmetric almost frame-like wavelet system providing approximation order $n$.
	A transformation based on the lifting scheme  is described that allows to improve almost frame-like wavelets to dual wavelet frames and preserve other properties. 
	A direct method for the construction of dual wavelet frames providing approximation order $n$ and mutual symmetry properties is also discussed.
	For an abelian symmetry group ${\cal H}$,
	a technique providing the ${\cal H}$-symmetry property for each wavelet function is given for the above three methods. 
\end{abstract}


\end{frontmatter}


\section*{Introduction}

Multivariate nonseparable wavelet systems have its applications in 2-D tomography, 3-D rotational angiography and other fields (see, e.g.\cite{Gan},\cite{KST},~\cite{KAT}, \cite{RRRR} and the references therein). Such wavelets are more suitable for  applications than separable wavelets.
But in the multivariate case, the problem of the construction of wavelets with desired properties is more complicated comparing to the univariate case. 
One of the main difficulties is that there is no simple and effective algorithm for the matrix extension problem.
A number of papers are devoted to the construction of nonseparable wavelets for some concrete situations.

A symmetry property for wavelets is highly desired in applications since symmetric wavelets efficiently work with the edges of images and reduce the amount of computations. Symmetric wavelets are connected with linear-phase filter banks, which are more applicable for images since they produce no phase distortion.
But due to the different kinds of symmetry,
it is more complicated to construct symmetric wavelets in general case,
in contrast to the univariate case where the general techniques are known (see~\cite{HanZh13},~\cite{Pet}).
In the multivariate case methods for the
construction of point symmetric wavelets can be found in~\cite{GLS} and~\cite{KST}.
Some schemes in different setups for the construction of highly symmetric wavelet systems were presented in~\cite{And},~\cite{Jiang4fold},~\cite{Jiang6fold},~\cite{Koch} and the references therein.
The construction of wavelet masks  for the interpolatory case was considered in~\cite{KrSymInterp}.

The paper is devoted to the construction of 
multivariate wavelet frames with
important for applications properties, i.e. good approximation order and symmetry properties. Let ${\cal H}$ be a symmetry group, $M$ be a dilation matrix.
For a given ${\cal H}$-symmetric refinable mask obeying the sum rule of order $n$, in Theorem~\ref{theoWaveNoSym} an algorithm  is suggested for the
construction of  almost frame-like wavelet system providing approximation order $n$ and having mutual symmetry properties.
The number of wavelet generators is equal to $m$ or $m-1$, $m=|\det M|.$ 
Frame-like wavelets were introduced in~\cite{KrSk}. These systems are not dual wavelet frames but preserve important properties of frames. Among the benefits, there is a simplification of the construction, since we do not need to provide vanishing moments for all wavelet functions (this condition is required for dual wavelet frames).
Using the lifting scheme it is possible to improve almost frame-like wavelet system to dual wavelet frame, preserving the symmetry and approximation properties. 
Also, we discuss a direct approach for the construction of dual wavelet frames.
This approach  is based on Algorithm 1 from~\cite{Sk1}. 
In this case, for a given ${\cal H}$-symmetric refinable mask obeying the sum rule of order $n$, in Theorem~\ref{theoFrameNoSym} the method is given for the construction of dual wavelet frames such that  wavelet masks are mutually symmetric and have vanishing moments of order $n$. The number of wavelet generators, in this case, is equal to $m$ or $m+1.$ But the support of the constructed wavelet functions can be large and this is not good for applications.
For an abelian symmetry group ${\cal H}$ and under some additional assumptions, a symmetrization step can be done for the three setups described above. As a result, all wavelet masks have the ${\cal H}$-symmetry properties keeping other properties unchanged. 

The paper is organized as follows.  In Section 1 we give some basic notations and definitions including the notion of a symmetry group and connected notions.
Section 2 is devoted to  the construction of wavelet masks with mutual symmetry properties for three cases: frame-like wavelets, frames based on the lifting scheme, frames based on Algorithm 1 from~\cite{Sk1}.
Section 3 describes a symmetrization step which allows to construct fully  ${\mathcal H}$-symmetric wavelet systems in the same three cases.
In Section 4 several examples are presented.

\section{Basic notations and definitions}

	We use the standard multi-index notations.
	For $x,y\in\rd,$
	$(x, y)=\sum\limits_{i=1}^d x_iy_i$,
	${\bf0}=(0,\dots, 0)\in{\mathbb R}^d$
	and
	$x \geq y$ if $x_j \geq y_j,$ $j= 1,\dots, d.$
	${\mathbb Z}_+^d:=\{x\in{\mathbb Z}^{d}:~x\geq~{\bf0}\}.$
	If $\alpha,\beta\in{\mathbb Z}^{d}_+$, $b\in{\mathbb R}^d$, we set
	$|\alpha|=\sum\limits_{j=1}^d \alpha_j$,
	$\alpha!=\prod\limits_{j=1}^d\alpha_j!$,
	$\binom{\alpha}{\beta}=\frac{\alpha!}{\beta!(\alpha-\beta)!}$,
	$b^{\alpha}=\prod\limits_{j=1}^d b_j^{\alpha_j}$,
	$D^{\alpha}f=\frac{\partial^{|\alpha|} f}{\partial x^{\alpha}}.$
	For $n\in{\mathbb N}$, $\Delta_n~:=~\{\alpha\in{\mathbb Z}^{d}_+:\,\,|\alpha|<n\}.$
	For $\beta\in\zd_+$, $\square_{\beta}:=\{\alpha\in\zd_+:\,\alpha\le\beta\}.$
	For a set $X,$  $\# X$ denotes the cardinality of $X.$
	For a matrix $M$, $[M]_{i,j}$ denotes 
	the (i,j)-th element of $M$,
	$[M]^i$ denotes the i-th column of $M$,
	$M^*:=\overline{M^T}$.
	$I_d$ denotes the $d\times d$ identity matrix.

	An  integer  $d\times d$ matrix $M$ is called a\textit{ dilation matrix} if
	its eigenvalues  are greater than 1 in modulus, $m=|\det M|.$ 
	A complete set of representatives  of the quotient group $\zd \slash M \zd$ is denoted by $D(M)=\{s_0,\dots, s_{m-1}\}$. An element of $D(M)$ is called a \textit{digit},
	$D(M)$ is called \textit{the set of digits}
	(see, e.g., \cite[\S~2.2]{NPS}). We assume that $s_0={\bf0}.$

	The Fourier transform of $f\in L_1({\mathbb R}^d)$
	is 
	$\widehat
	f(\xi)=\int_{{\mathbb R}^d} f(x)e^{-2\pi i (x, \xi)}\,dx,$ $\xi\in{\mathbb R}^d.$
	This notion can be naturally extended to $L_2({\mathbb R}^d)$
	and     to the space of tempered distributions ${\cal S}'.$
	For a function $f$ defined on $\rd$ set $f_{jk} : = m^{j/2} f(M^j \cdot + k),$
	$j\in\z$, $k\in\zd.$	
	
A finite set ${\mathcal H}$ of $d\times d$ unimodular matrices (i.e. integer matrices with determinant equal to $\pm 1$) is 
	a \emph{symmetry group}, if ${\mathcal H}$ forms a group under the matrix multiplication. We say that a dilation matrix $M$ is \textit{appropriate} for ${\mathcal H}$, if $M^{-1}EM\in {\mathcal H},$ $\forall E\in {\mathcal H}.$   Or, equivalently, for each $E\in {\mathcal H}$ there exists
	$E'\in {\mathcal H}$ such that
	\begin{equation}
	EM=ME' \quad \texttt{or} \quad M^{-1}E=E'M^{-1}.
	\label{fMEEM}
	\end{equation} 	
	A compactly supported distribution $f$ is called \emph{${\mathcal H}$-symmetric} with respect to  a
	center $C\in {\mathbb R}^d$, if
	\[
	\h f(\xi) = \h f(E^*\xi) e^{2\pi i (EC-C,\xi)},
	\quad\forall E\in {\mathcal H}, \quad \xi\in{\mathbb R}^d.
	\]
	For trigonometric polynomials we use a bit different definition.
	We say that  $c\in{\mathbb R}^d$ is an \textit{appropriate symmetry center} for ${\mathcal H}$, if  $c-Ec\in{\mathbb Z}^{d},$ $\forall E\in {\mathcal H}$.
	A trigonometric polynomial
	$t(\xi)=\sum_{k\in{\mathbb Z}^{d}}h_k e^{2\pi i (k, \xi)},$ $h_k \in{\mathbb C},$
	is  \emph{${\mathcal H}$-symmetric }
	with respect to an appropriate center $c$, if
	\be
	t(\xi)=e^{2\pi i (c-Ec,\xi)}t(E^*\xi), \quad
	\forall E\in {\mathcal H}.
	\label{TrigSym}
	\ee

	A function/distribution $\phi$ is called \emph{refinable} 
	if there exists a $1$-periodic
	function $m_0\in L_2([0,1]^d)$ 
	(mask, also refinable mask, low-pass filter) such that
	\be
	\widehat\phi(\xi)=m_0(M^{*-1}\xi) \widehat\phi(M^{*-1}\xi).	
	\label{RE}
	\ee
	This condition is called the refinement equation.
	It is well
	known (see, e.g.,~\cite[\S~2.4]{NPS}) 
	that for any  trigonometric polynomial
	$m_0$ satisfying $m_0({\bf0})=1$ there exists a unique 
	(up to a factor) compactly supported solution  of
	the refinement equation~(\ref{RE}) in the space of tempered distributions $S'$. 
	Throughout the paper we assume that any 
	refinable mask $m_0$ is a trigonometric
	polynomial  and $m_0({\bf0})=1$.
For an appropriate $M$ for ${\cal H}$, it is known that refinable mask $m_0$ is ${\cal H}$-symmetric with respect to an appropriate center $c$ if and only if the
corresponding refinable function is ${\cal H}$-symmetric with respect to $(M-I_d)^{-1}c$ (see~\cite{HanSym02}).

	Let us fix the set of digits $D(M)$. For any trigonometric polynomial $t$ there exists a unique set of
	trigonometric polynomials $\tau_k$, $k=0\ddd m-1$, such that
	\be
	t(\xi)= \frac1{\sqrt m}\sum\limits_{k=0}^{m-1} e^{2\pi i(s_k,\xi)}\tau_k(M^*\xi),
	\label{PR}
	\ee
	where $s_k \in D(M)$.
	Equality~(\ref{PR}) is the \emph{polyphase representation} of  $t$.
	Trigonometric polynomial $\tau_k$ is called the \emph{polyphase component} of $t$
	corresponding to the digit $s_k.$

	Let $t$ be a   trigonometric polynomial, $n\in\n.$
	We say that  $t$ obey the \emph{sum rule of order $n$}
	with respect to dilation matrix $M$, if
	$
	D^\beta t({M^*}^{-1}\xi)\Big|_{\xi=s}=0,$ $\forall s\in D(M^*)\setminus\{\nul\},$
	$\forall  \beta\in\Delta_n.
	$
	We say that  $t$ has \emph{vanishing moments of order $n$} if
	$D^{\beta} t(\nul)= 0,$   $ \forall  \beta\in\Delta_n.$

\subsection{Symmetry groups and digits}

	A general method for the construction of multivariate ${\cal H}$-symmetric 
	masks that obey the sum rule of an arbitrary order is given in~\cite{KrSym}. The key feature used in the construction is the description of how a symmetry group	acts on the set of digits. 
	Let a center $c$ and a dilation matrix $M$ be appropriate for ${\cal H}$.
	It is known that any $\alpha\in{\mathbb Z}^{d}$ can be 
	uniquely represented as $\alpha=M\beta+s,$
	where $\beta\in {\mathbb Z}^{d},$ $s\in D(M).$     
	This fact yields that for each digit $s\in D(M)$
	and matrix $E\in {\mathcal H}$ there exist a unique digit $q\in D(M)$ 
	and a unique vector $r_s^E\in {\mathbb Z}^{d}$ such that 
	\be
	Es=M r_{s}^E+q+Ec-c.
	\label{fDigMain}
	\ee 
	The indices of $r_s^E$ mean that vector $r_s^E$ depends on digit $s$ and matrix $E.$
	
	The coset corresponding to digit $s\in D(M)$ we denote by
	$\langle s\rangle$, i.e. 
	$\langle s\rangle=M\zd+s.$      Denote by ${\cal D}:=\{\langle s\rangle, s\in D(M)\}$ the set of cosets.    
	Define a group action of 
	${\cal H}$ on ${\cal D}$
	as follows 
	$$
	E\langle s\rangle:=\{EM\beta+Es+c-Ec, \beta\in{\mathbb Z}^{d}\}, \quad E\in {\cal H}, \quad s\in D(M).
	$$
	Note that $E\langle s\rangle$ is also a coset, i.e. there exists $q\in D(M)$ such that $E\langle s\rangle = \langle q\rangle$. 
	Next, we introduce suitable notations that will be used throughout the paper. These notations are illustrated in~\cite{KrSym} and~\cite{KrSymInterp}.
	\begin{itemize}
		\item  ${\cal H}\langle s \rangle=\{E \langle s \rangle, E\in {\cal H}\}$ is the \textit{orbit} of $\langle s \rangle\in {\cal D}$.
		The orbits are disjoint,
		${\cal H}\langle s \rangle \subseteq {\cal D}.$ 
		
		\item The set $\Lambda\subset {\cal D}$ contains representatives 
		from each orbit. ${\cal D}= \cup_{\langle s \rangle\in \Lambda}{\cal H}\langle s \rangle.$ 
		
		\item For convenience, redenote the elements of the set $\Lambda$ by $\langle s_{p,0}\rangle,$ where $p=0,\dots,\#\Lambda-1.$ 
		
		\item The set
		$
		{\cal H}_{p,0}=\{F\in {\cal H} : F\langle s_{p,0} \rangle=\langle s_{p,0} \rangle\}
		$
		is the  \textit{stabilizer} of $\langle s_{p,0}\rangle $; 
		${\cal H}_{p,0} \subset {\cal H}.$
		
		\item The set ${\cal E}_p$ contains a complete set of representatives of ${\cal H}\slash {\cal H}_{p,0};$ 
		${\cal E}_p \subset {\cal H}.$
		
		\item The elements of the orbit ${\cal H}\langle s_{p,0}\rangle$ we denote by $\langle s_{p,i}\rangle,$ $i=0,\dots,\#{\cal E}_p-1.$
		
		\item For a fixed index $p$, the matrices of the set 
		${\cal E}_p$ we denote by $E^{(i)}$ 
		such that $E^{(i)} \langle s_{p,0}\rangle= \langle s_{p,i}\rangle, $ $i=0,\dots,\#{\cal E}_p-1.$ Note that $E^{(0)}=I_d.$
		
		\item The digit corresponding to the coset $\langle s_{p,i}\rangle$ we choose such that
		\be
		E^{(i)} s_{p,0}+c-E^{(i)} c=:s_{p,i}, 
		\quad i=1,\dots,\#{\cal E}_p-1.
		\label{fDigE_i}
		\ee
	\end{itemize}

	Note that for a fixed $p$,
	$p=0,\dots,\#\Lambda-1$, 
	symmetry group ${\cal H}$ can be uniquely represented as
	follows ${\cal H}={\cal E}_p \times {\cal H}_{p,0},$ i.e. for each matrix $K$ in ${\cal H}$ there exist matrices $E\in {\cal E}_p$ and $F\in {\cal H}_{p,0}$ such that $K=E F.$  
	The sets ${\cal E}_p, $ ${\cal H}_{p,0}$ can be considered as the\textit{ "coordinate axes"} of  symmetry group ${\cal H}$. 
	For each $p$ 
	these "coordinate axes" of ${\cal H}$ can be different.
	Also note how ${\cal H}$ acts on a digit $s_{p,0}\in D(M)$.
	If $E^{(i)}\in{\cal E}_p$, then  $E^{(i)} s_{p,0}$ is defined in~(\ref{fDigE_i}).
	If $F\in {\cal H}_{p,0}$, then $F\langle s_{p,0}\rangle=\langle s_{p,0}\rangle$ and
	\be
	Fs_{p,0}=M r_{p,0}^F+ s_{p,0}+Fc-c,
	\label{fDigF}
	\ee
	where $r_{p,0}^F\in {\mathbb Z}^{d}$. 
	Notice that $ r_{p,0}^F= M^{-1}(c-s_{p,0})-M^{-1}F(c-s_{p,0}).$
	
	%
	%

	Now we formulate the 
	${\cal H}$-symmetry condition~(\ref{TrigSym}) for a trigonometric polynomial in terms of its polyphase components.
Alongside with the standard enumeration of the polyphase components of trigonometric polynomial $t$, we also use enumeration corresponding to the new enumeration of digits: $\tau_{p,i}(\xi),$
$i=0,\dots,\#{\cal E}_{p}-1,$
$p=0,\dots,\#\Lambda-1$.

	\begin{lem}~\cite[Lemma 8]{KrSymInterp}
		A trigonometric polynomial $t$ is ${\cal H}$-symmetric with respect to
		an appropriate center $c$
		if and only if
		its polyphase components $\tau_{p,i}$  
		satisfy 
		\be
		\tau_{p,0}(\xi)=
		e^{2\pi i (r^F_{p,0},\xi)}\tau_{p,0}((M^{-1} F M)^*\xi),
		\qquad  \mbox{for all } F\in {\cal H}_{p,0};
		\label{fPolyF}
		\ee
		\be          
		\tau_{p,i}(\xi)=\tau_{p,0}((M^{-1} E^{(i)} M)^*\xi),
		\qquad  E^{(i)}\in{\cal E}_{p},
		\quad    \mbox{for all }     i\in\{0,\dots,\#{\cal E}_{p}-1\}.
		\label{fPolyE_i}
		\ee
		for each $p\in\{0,\dots,\#\Lambda-1\}$.
		\label{lemSymPoly}
	\end{lem}

%

\section{Construction of symmetric wavelets}

	A family of functions $\{f_\alpha\}_{\alpha\in\aleph}$ ($\aleph$ is a countable index set) in a
	Hilbert space $H$ is called a frame in $H$
	if there exist constants $A, B > 0$ such that  
	$A\|f\|^2\le\sum\limits_{\alpha\in\aleph}|\langle f, f_\alpha\rangle|^2\le B \|f\|^2,$ $\forall f\in H. $
	If $\{f_\alpha\}_\alpha$ is a frame in $H$, then
	every $f\in H$ can be decomposed as
	$ f=\sum_\alpha\langle f,\w f_\alpha\rangle f_\alpha, $ where $\{\w f_\alpha\}_\alpha$ is a
	dual frame in $H$. Comprehensive characterization of frames
	can be found in~\cite{Chris}.  Wavelet frames are of great interest in many applications, especially in signal processing. For more information about multivariate wavelet frames see~\cite{SKP}.

For $\psi^{(\nu)}\in S'$, $\nu=1\ddd r$, a system
$\{\psi_{jk}^{(\nu)}\}$ is called a wavelet system.
We say that wavelet system
$\{\psi_{jk}^{(\nu)}\}$ has
\emph{vanishing moments of order $n \in\n$}
(or has the {\em  $VM^{n}$ property})
if
$D^{\beta}\h{\psi^{(\nu)}}({\bf0})=0,$  $\forall \beta\in\Delta_n,$
$\nu=1\ddd r,$

A general scheme for the construction of compactly supported MRA-based wavelet systems (in particular, wavelet frames in $L_2(\rd)$) 
was developed in~\cite{RS},
\cite{RS2} (the Unitary Extension Principle).
To construct a pair of such wavelet systems  
one
starts with two compactly supported refinable functions
$\phi, \widetilde\phi$, $\h\phi(\nul)=1$, $\h{\w\phi}(\nul)=1,$
(or its masks $m_0, \widetilde m_0$, respectively, which are trigonometric
polynomials).
Then one finds trigonometric
polynomials $m_{\nu}, \widetilde m_{\nu}$, $\nu = 1\ddd r$,
$r\ge m-1$, called \textit{wavelet masks}, such that the following polyphase matrices 
$$
{\cal M}:=\{\mu_{\nu k}\}_{\nu={0,r}}^{k={0,m-1}},\,
\widetilde{\cal  M}:=\{\w\mu_{\nu k}\}_{\nu={0,r}}^{k={0,m-1}}
\quad \mbox{ satisfy } \quad     
{\cal M}^*{\cal\widetilde M}= I_m.
$$
Here $\mu_{\nu k},$ $\w\mu_{\nu k}$ $k=0,\dots,m-1,$ are the polyphase components of the
wavelet masks $m_{\nu},$ $\w m_{\nu}$ for all $\nu=0,\dots,r,$ $r\ge m-1.$
The wavelet functions $\psi^{(\nu)}$, $\w\psi^{(\nu)},$
$\nu = 1\ddd r$, $r\ge m-1,$ are defined via its Fourier transform
$ \widehat{\psi^{(\nu)}}(\xi)=
m_\nu(M^{*-1}\xi)\widehat\varphi(M^{*-1}\xi),$ 
$\widehat{\widetilde\psi^{(\nu)}}(\xi)=
\widetilde m_\nu(M^{*-1}\xi)\widehat{\widetilde\varphi}(M^{*-1}\xi).
$
If the wavelet functions  $\psi^{(\nu)}, \widetilde\psi^{(\nu)}$,
$\nu=1\ddd r$, $r\ge m-1,$
are constructed as above, then
the set of the functions $\{\psi_{jk}^{(\nu)}\}$, $\{\w\psi_{jk}^{(\nu)}\}$
is said to be a
\emph{compactly supported MRA-based dual wavelet system} generated by the
refinable functions $\varphi, \widetilde\varphi$
(or their masks $m_0, \widetilde m_0$).

Suppose that some  compactly supported refinable functions $\varphi, \widetilde\varphi\in L_2(\rd)$ generate a MRA-based dual wavelet system $\{\psi_{jk}^{(\nu)}\}$, $\{\w\psi_{jk}^{(\nu)}\},$
$\nu=1\ddd r$, $r\ge m-1.$
A necessary (see~\cite[Theorem 1]{Sk1}) and
sufficient (see~\cite[Theorems 2.2, 2.3]{Han1}) condition
for the wavelet system $\{\psi_{jk}^{(\nu)}\}$, $\{\w\psi_{jk}^{(\nu)}\}$
to be a pair of dual wavelet frames in $L_2(\rd)$ is that
each wavelet system $\{\psi_{jk}^{(\nu)}\}$ and $\{\w\psi_{jk}^{(\nu)}\}$
has vanishing moments at least of order $1$.
Or, equivalently, the corresponding wavelet masks
$m_{\nu},$ $\w m_{\nu}$
have vanishing moments at least of order $1$.
Also, it is known that good approximation properties 
for the corresponding wavelet systems are provided
by the $VM^{n}$ property for the dual wavelet system $\{\w\psi_{jk}^{(\nu)}\}$ (see, e.g.,~\cite[Theorem 4]{Sk1}).
Although algorithms for the construction of dual wavelet frames
with vanishing moments were developed (see,~\cite{Sk1}), 
but it is not easy to provide various types of symmetry for different  dilation matrices.    
If we reject the frame requirements and aim to provide the 
$VM^{n}$ property only for dual wavelet system $\{\w\psi_{jk}^{(\nu)}\}$, $\nu=1\ddd r$, then
the method can be simplified. The necessary conditions with the constructive proof are given in

\begin{lem} \cite[Lemma 14]{KrSk} Let  $\phi, \widetilde\phi$ be compactly supported  refinable distributions,
	$m_0$ obeys the sum rule of order $n$ and
	\be
	D^\beta\lll 1-m_0(\xi)\overline{\w m_0(\xi)}\rrr\Big|_{\xi=0}=0 \qquad
	\forall \beta\in\Delta_n.
	\label{20_new}
	\ee
	then there exist a MRA-based dual wavelet system
	$\{\psi_{jk}^{(\nu)}\}$, $\{\w\psi_{jk}^{(\nu)}\}$, $\nu=1\ddd m,$
	such that  wavelet system $\{\widetilde\psi^{(\nu)}_{jk}\}$
	has the $VM^{n}$ property.
	\label{lemKrSk}
\end{lem}

A technique for the extension of the polyphase matrices
realizes as follows. Let us denote the rows of the polyphase components of masks $m_0$ and $\w m_0$ by $P=(\mu_{00},\dots,\mu_{0,m-1}),$ 
$\w P=(\widetilde\mu_{00},\dots,\widetilde\mu_{0,m-1})$ and extend them with the elements $1-P\w P^*$ and $1$ accordingly. Then  the explicit formulas  for the  matrix extension are given by 
\be
{\cal N}=\left( \begin{array}{cc}
	P & 1-P\w P^*\\
	U & -U\w P^*
\end{array}\right),
\quad {\cal \w N}=\left(
\begin{array}{cc}
	\w P& 1\\
	\w U-\w U P^* \w P& -\w U P^*
\end{array} \right).
\label{fTTblock}
\ee
Here $U,$ $\w U$ are $m\times m$ matrices consisting of trigonometric polynomials such that
$U\w U^*\equiv I_m.$ For instance, we can take $U=\w U=I_m.$ It follows easily that  ${\cal N}\w {\cal N}^*\equiv I_{m+1}.$ 
It is important to note that dual wavelet frames cannot be constructed using such extension technique due to
the fact that in the matrix $\w{\cal N}$ 
the upper right element is equal to one. 
Therefore,  wavelet masks $m_{\nu},$ $\nu=1,\dots,m$ do not have any order of vanishing moments
according with~\cite[Theorem 8]{Sk1}.
Nevertheless, for appropriate functions $f$ we can consider frame-type decompositions    with respect to this system.

Let $\{\psi_{jk}^{(\nu)}\}$, $\{\w\psi_{jk}^{(\nu)}\}$  be a MRA-based dual
wavelet system and  $A$ be a class of functions $f$ such that
$\langle f,\widetilde\phi_{0k}\rangle$,
$\langle f,\widetilde\psi_{jk}^{(\nu)}\rangle$ have meaning 
(e.g.,  $f$ in  the Schwartz space $S$ if $\phi,\,\widetilde\psi^{(\nu)}\in S'\,$;
$f\in L_p\,$ if $\phi, \widetilde\psi^{(\nu)}\in L_q$, $\frac1p+\frac1q=1$).
We say that
dual wavelet system $\{\psi_{jk}^{(\nu)}\}$ $\{\w \psi_{jk}^{(\nu)}\}$  is {\em almost frame-like} if
\be
f=\sum\limits_{k\in\,\zd}\langle f,\widetilde\phi_{0k}\rangle\phi_{0k}+
\sum\limits_{j=0}^{\infty}\sum\limits_{\nu=1}^r\sum\limits_{k\in\,\zd}
\langle f,\widetilde\psi_{jk}^{(\nu)}\rangle\psi_{jk}^{(\nu)}\quad \forall f\in A,
\label{02}
\ee
where the series in (\ref{02}) converge in some natural sense.
Here we give the result from~\cite{KrSk}.

\begin{theo} \cite[Theorem 12]{KrSk} Let $f\in S$, $\phi,\widetilde\phi\in S'$,
	$\phi, \widetilde\phi$ be compactly supported and refinable,
	$\h\phi(\nul)=\h{\widetilde\phi}(\nul)=1$.
	Then the MRA-based dual wavelet system $\{\psi_{jk}^{(\nu)}\}$, $\{\w\psi_{jk}^{(\nu)}\},$
	$\nu=1\ddd r,$
	generated by $\phi, \widetilde\phi$
	is almost frame-like, i.e.~(\ref{02}) holds
	with the series converging in  $S'$.
	\label{t6}
\end{theo}

If, moreover, $\phi, \widetilde\phi$ are as in Lemma~\ref{lemKrSk}, $\phi\in L_2(\rd)$, then
(\ref{02}) holds with the series converging in $L_2$-norm and the corresponding almost frame-like system has approximation order $n$ (see \cite[Theorem 16]{KrSk}).
In the next sections we show that starting from symmetric refinable masks,  
the symmetry property for the wavelet masks can be provided.  
%

\subsection{Dual masks}

	Given an ${\cal H}$-symmetric refinable mask $m_0$,
	the first step for the construction of symmetric wavelets is the construction of an appropriate ${\cal H}$-symmetric dual refinable mask $\w m_0$.
	The method is given by the following theorem.
	
	
	\begin{theo}
		Let a dilation matrix $M$ and a center $c$
		be appropriate for a symmetry group ${\cal H},$ 
		$n\in\n$.
		Suppose $m_0$ is an ${\cal H}$-symmetric
		with respect to the center $c$ mask that obeys the sum rule of order $n$.
		Then there exists a dual mask $\w m_0$ that is ${\cal H}$-symmetric
		with respect to the center $c$ and satisfies condition~(\ref{20_new}).
		\label{DualTheo}
	\end{theo}
	
	{\bf Proof.} 
	Set $(2\pi i)^{|\beta|} \lambda_{\beta}:=D^\beta  m_0(\nul)$,
	$\beta\in\Delta_n$. And let numbers $\w \lambda_\beta$ satisfy
	\be
	\sum_{\alpha\in\square_\beta} 
	(-1)^{|\beta-\alpha|}        
	\binom{\beta}{\alpha}
	\lambda_{\alpha}\overline{ \w \lambda_{\beta-\alpha}} = 0, \quad
	\forall \beta \in \Delta_n.
	\label{fLambda}
	\ee
	Numbers  $\w \lambda_\beta$ can be found recursively from~(\ref{fLambda}).
	Define mask $\w m_0$  as follows
	\be
	\w m_0(\xi)=
	\frac 1 {\# {\cal H}}
	\sum\limits_{E\in {\cal H}} G(E^*\xi)e^{2\pi i (c-Ec,\xi)},
	\label{fDualMaskA}
	\ee
	where    $G(\xi)$  is a trigonometric polynomial such that 
	$ D^\beta G(\nul)=
	(2\pi i)^{|\beta|}\w\lambda_{\beta}$ for all $\beta\in\Delta_n$.
	It is not hard to check that $\w m_0(\xi) $ is ${\cal H}$-symmetric with respect to the point $c$. 
	Let us show that condition~(\ref{20_new}) is also valid.
	When $\beta=\nul$, then $\lambda_\nul = \w\lambda_\nul = 1,$
	and $m_0 (\nul) = \w m_0(\nul) = 1.$
	For $\beta \neq \nul$ condition~(\ref{20_new}) will be valid if
	\be  
	D^\beta\lll m_0(\xi)
	\overline{ G(E^*\xi)e^{2\pi i (c-Ec,\xi)}}\rrr\Big|_{\xi=0}=0,
	\quad \forall \beta\in\Delta_n, \beta \neq \nul, \quad \forall
	E\in {\cal H}.
	\label{fM0WM0}
	\ee  
	Since $m_0$ is ${\cal H}$-symmetric, then $m_0(\xi)
	\overline{ G(E^*\xi)e^{2\pi i (c-Ec,\xi)}} = m_0(E^*\xi)
	\overline{ G(E^*\xi)}$. Due to the higher chain rule with
	the linear change of variables, equalities~(\ref{fM0WM0}) are equivalent to
	$$
	D^\beta\lll m_0(\xi)
	\overline{ G(\xi)}\rrr\Big|_{\xi=0}=0,
	\quad \forall \beta\in\Delta_n, \beta \neq \nul.
	$$
	These equalities are valid, since $D^\beta G(\nul)=
	(2\pi i)^{|\beta|} \w\lambda_{\beta}$,  $\forall\beta\in\Delta_n$, and numbers $\w\lambda_{\beta}$ 
	satisfy~(\ref{fLambda}).
	$\Diamond$

	For applications it is important to get  refinable  masks 
	with the minimal coefficient support 
	(the coefficient support of trigonometric polynomial $t=\sum_k h_k e^{2\pi i (k,\xi)}$ is defined by   $\texttt{coefsupp} (t) = \{k\in\zd: h_k\neq 0\}$).
	The general recipe for the construction of an ${\cal H}$-symmetric dual refinable mask with the minimal coefficient support  is based on a fact that
	for any vector $k\in \texttt{coefsupp} (G)$, the set of integer   points
	$\{E k + c - Ec : E\in {\mathcal H}\}$ should be a subset of $ \texttt{coefsupp} (\w m_{0})$ (see~(\ref{fDualMaskA})).
	Thus, for some fixed $k_1\in\zd$, one should check whether the system of linear equations $D^\beta G(\nul)=
	(2\pi i)^{|\beta|} \w\lambda_{\beta}$,  $\forall\beta\in\Delta_n$, is solvable with
	$\texttt{coefsupp} (G) =K_1:= \{E k_1 + c - Ec: E\in {\mathcal H}\}$.
	If yes, then one solves the system. If not, one should add another point
	$k_2\in\zd$ and check the solvability of the system with 
	$\texttt{coefsupp} (G) =K_1 \cup \{E k_2 + c - Ec: F\in {\mathcal H}\}$ and so on.

\subsection{Symmetric frame-like wavelets}
\label{sec:SymFL}

	Assume we have two ${\cal H}$-symmetric masks $m_0$ and $\w m_0$ such that~(\ref{20_new}) is valid.
	Let us consider the construction of wavelet systems using the basic matrix extension~(\ref{fTTblock}) with $U=\w U=I_{m}.$ Namely,
	\be
	{\cal N}=\left( \begin{array}{c|c}
		P & 1-P\w P^* \\
		\hline
		&  \\
		I_m & -\w P^{*}
	\end{array}\right),
	\quad 
	{\cal \w N}=\left(
	\begin{array}{c|c}
		\w P& 1 \\
		\hline
		&   \\
		I_m- P^{*} \w P & -P^{*}
	\end{array} \right),
	\label{fNNblockNoSym}
	\ee
	where $P$ and $\w P$ are the rows of the polyphase components of masks $m_0$ and $\w m_0.$
	For convenience, we use the following enumeration of wavelet masks.
	Let $\mu_{0,p,i}$ and $\w\mu_{0,p,i}$, $i=0,\dots,\#{\cal E}_{p}-1,$     $p=0,\dots,\#\Lambda-1,$ be the polyphase components of $m_0$ and $\w m_0$,
	$$T_p=(\mu_{0,p,0},...,\mu_{0,p,\#{\cal E}_{p}-1}), \qquad
	\w T_p=(\w\mu_{0,p,0},...,\w\mu_{0,p,\#{\cal E}_{p}-1}).$$
	Set $P=(T_0,\dots,T_{\#\Lambda-1}),$ 
	$\w P=(\w T_0,\dots,\w T_{\#\Lambda-1})$.
	The polyphase components for wavelet masks $m_\nu$ and 
	$\w m_\nu$ are contained in the submatrices $I_m$  and $I_m- P^{*} \w P$ in~(\ref{fNNblockNoSym}).
	Note that $I_m- P^{*} \w P$ is a block matrix 
	\be
	I_m-P^{*} \w P=
	\left( \begin{array}{cccc}
		I_{N_0-1}-
		T_0^* \w T_0 & 
		-T_0^*\w T_1 & \dots & 
		-T_0^*\w T_{\#\Lambda-1}\\
		-T_1^*\w T_0 &
		I_{N_1-1}-T_1^*\w T_1 & \dots & 
		-T_1^*\w T_{\#\Lambda-1} \\
		\vdots & & \ddots & \vdots \\
		-T_{\#\Lambda-1} \w T_0 &	\dots &\dots &
		I_{N_{(\#\Lambda-1)}-1}-T_{\#\Lambda-1} 
		\w T_{\#\Lambda-1}
	\end{array}\right)
	\label{eq:Nblock}    
	\ee
	where $N_p = \# {\cal E}_p$.
	The rows of the submatrix we enumerate by double index $(p,i)$.
	The  $(p,i)$th row is the $i$th row in the $p$th block.
	A wavelet mask corresponding to the $(p,i)$th  row is denoted by
	$\w m_{(p,i)}$.
	Analogously, wavelet masks $m_{(p,i)}$ are enumerated.

	\begin{theo} 
		Let a dilation matrix $M$ and a center
		$c$ be  appropriate for a symmetry group ${\cal H}$, $n\in\n$.
		Suppose $m_0$ and $\w m_0$ are ${\cal H}$-symmetric
		with respect to the center $c$ masks such that
		mask $m_0$
		obeys the sum rule of order $n$,
		mask $\w m_0$ satisfies
		condition~(\ref{20_new}).
		Then wavelet masks  $m_{(p,i)}$ and $\w m_{(p,i)}$
		constructed using matrix extension~(\ref{fNNblockNoSym})
		have the following symmetry properties:
		\begin{enumerate}[{\bf W\arabic*}]
			\item \label{fLS_1cond} $m_{(p,0)}$, $\w m_{(p,0)}$ are ${\cal H}_{p,0}$-symmetric 
			with respect to   the center $s_{p,0}$,
			\item \label{fLS_2cond}  $m_{(p,i)}(\xi)=m_{(p,0)}(E^{(i)*}\xi) e^{2\pi i (c-E^{(i)}c,\xi)},$
			
			$\w m_{(p,i)}(\xi)=\w m_{(p,0)}(E^{(i)*}\xi) e^{2\pi i (c-E^{(i)}c,\xi)},$  
			$E^{(i)}\in{\cal E}_{p},$
		\end{enumerate}
		where  $i=0,\dots,\#{\cal E}_{p}-1,$
		$p=0,\dots,\#\Lambda-1$.
		Wavelet masks $\w m_{(p,i)}$ have vanishing moments of order $n$. 
		The corresponding MRA-based dual wavelet system 
		is almost frame-like in $S'.$    
		\label{theoWaveNoSym}
	\end{theo}
	
	{\bf Proof.} Wavelet masks from submatrix $I_m- P^{*} \w P$ in~(\ref{fNNblockNoSym}) are given by
	$$\w m_{(p,i)}(\xi)=\frac{1}{\sqrt{m}}e^{2\pi i (s_{p,i},\xi)}-\overline{\mu_{0,p,i}(M^*\xi)}\w m_0(\xi), \quad  i=0,\dots,\#{\cal E}_{p}-1,
	p=0,\dots,\#\Lambda-1.$$
	Then mask $m_{(p,0)}$ is ${\cal H}_{p,0}$-symmetric with respect to the center $s_{p,0}.$
	Indeed, due to~(\ref{fDigF}) and~(\ref{fPolyF})
	for all  $F\in {\cal H}_{p,0}$  we have
	\begin{equation*}
	\begin{split}
		\w m_{(p,0)}(F^*\xi)  = &
		\frac{1}{\sqrt{m}} e^{2\pi i (F s_{p,0},\xi)} - \overline{\mu_{0,p,0}(M^*F^*\xi)} \w m_0(F^*\xi) \\ =&
		\w m_{(p,0)}(\xi)e^{2\pi i (Fs_{p,0}-s_{p,0},\xi)}.
	\end{split}
	\end{equation*} 

	Also, by~(\ref{fDigE_i}) and~(\ref{fPolyE_i}) we have
		\begin{equation*}
	\begin{split} \w m_{(p,0)}(E^{(i)*}\xi) =  &
	\frac{1}{\sqrt{m}} e^{2\pi i (E^{(i)} s_{p,0},\xi)} - 
	\overline{\mu_{0,p,0}(M^*E^{(i)*}\xi)} \w m_0(E^{(i)*}\xi) \\ = &
	\w m_{(p,i)}(\xi)e^{2\pi i (E^{(i)}c-c,\xi)}.
	\end{split}
\end{equation*} 
	Similarly, it can be checked that the same symmetric properties are valid for  wavelet masks $ m_{(p,i)}$, since 
	$$ m_{(p,i)}(\xi)=\frac{1}{\sqrt{m}} e^{2\pi i (s_{p,i},\xi)}, \quad  i=0,\dots,\#{\cal E}_{p}-1,
	p=0,\dots,\#\Lambda-1.$$
	Vanishing moments of order $n$  for the dual wavelet masks $\w m_{(p,i)}$ are  provided by Lemma~\ref{lemKrSk}.
	$\Diamond$

	
	Thus, the matrix extension~(\ref{fNNblockNoSym}) leads to the wavelet masks that are\textit{ mutually symmetric}, i.e. some wavelet masks are reflected or rotated copies of the others. Note that~\ref{fLS_1cond} and~\ref{fLS_2cond} by direct computations imply that  
	$m_{(p,i)}$ and $\w m_{(p,i)}$ are ${\cal H}_{p,i}$--symmetric 
	with respect to   the center $s_{p,i}$. Proposition 2.1 in~\cite{Han3} allows to compute the symmetry centers of the corresponding wavelet functions. 
	Notice that if for some 
	${\cal H}$ and $M$ we get that ${\cal E}_p=\{I_d\}$ for all 
	$p=0,\dots,\#\Lambda-1$, then all wavelet masks constructed by
	Theorem~\ref{theoWaveNoSym}  are ${\cal H}$-symmetric.

The number of wavelet generators for  almost frame-like wavelet systems in Theorem~\ref{theoWaveNoSym} can be reduced to $m-1$ in the following case.
Suppose for some $k=0,\dots, m-1$, $\mu_{0k} \equiv \frac{1}{\sqrt m}$. Then  we can take $\w \mu_{0k} \equiv \sqrt m$ and $\w\mu_{0l} \equiv 0$, for $l\neq k$.  Therefore, $\w m_0(\xi) = e^{2\pi i (s_k,\xi)}$  and by Remark 15 in~\cite{KrSk}, since $\sum\limits_{l=0}^{m-1} \mu_{0l}(\xi) \overline{\w\mu_{0l}(\xi)} \equiv 1$, condition~(\ref{20_new}) is valid. The matrix extension~(\ref{fNNblockNoSym}) leads to a wavelet system with $m-1$ wavelet generators.

From the point of view of Theorem 1.1 in~\cite{HanShen}, the  wavelet system constructed by Theorem~\ref{theoWaveNoSym} can be a dual wavelet frame in a pair of dual Sobolev spaces depending on the Sobolev smoothness exponents of $\phi$ and $\w\phi$.

\subsection{Lifting scheme: frame-like wavelets to frames}

	The lifting scheme was introduced by W.~Sweldens~\cite{Sweld}.
	It is a tool for designing wavelets and performing  the discrete wavelet transform.
	The lifting scheme  has lots of different applications and useful properties.
	One important for us feature  is that the lifting scheme allows
	to improve properties of a given wavelet system. 
	In particular, it allows to provide additional 
	vanishing moments for wavelet masks (see, e.g.,~\cite{Bhatt}).
	Since vanishing moments are necessary and sufficient conditions
	for a wavelet system to be a dual wavelet frame, the lifting scheme
	can help to improve  frame-like wavelet systems
	to frames.
	In this section, we also give a method that helps to preserve symmetry properties
	of wavelets during the improvement.
	
	Let refinable masks $m_0$,  $\w m_0$ and  wavelet masks  $m_{\nu},$ $\w m_{\nu},$ $\nu=1,\dots,r$ be
	such that the corresponding polyphase matrices 
	${\cal M}$,
	$\widetilde{\cal  M}$
	satisfy $
	{\cal M}^*{\cal\widetilde M}= I_m.
	$
	Let $L_1,\dots,L_r$ be trigonometric polynomials.       
	Define $r\times r$ matrices ${\cal L},$ $\w {\cal L}$ as follows:
	$$
	{\cal L}:=\left(\begin{array}{cccc}
	1 & \nul\cr
	L^T & I_{r-1}  \cr        
	\end{array}\right),\,
	\w {\cal L}:=\left(\begin{array}{cccc}
	1 & -\overline{L}\cr
	\nul &  I_{r-1}           
	\end{array}\right),
	$$
	where $L=(L_1,\dots,L_r)$.
	It is easy to check that $
	{\cal L}^*{\cal\widetilde L}= I_r.
	$
	Define new polyphase matrices 
	${\cal M}_{new}:={\cal L}{\cal M}$ and
	$\w{\cal M}_{new}:=\w{\cal L} \w{\cal M}$.
	Note that the equality 
	$
	{\cal M}_{new}^*{\cal\widetilde M}_{new}= I_m
	$
	is preserved.
	Denote the elements of the new matrices as follows:
	${\cal M}_{new}=\{\mu^n_{\nu j}\}_{\nu=0,r}^{j=0,m-1}$,
	$\w{\cal M}_{new}=\{\w\mu^n_{\nu j}\}_{\nu=0,r}^{j=0,m-1}.$
	The transformation formulas for the polyphase components are the following:
	\begin{equation*}
	\begin{alignedat}{2}
	\mu_{0,j}^n(\xi)&=\mu_{0,j}(\xi), 
	& \quad &\mu_{\nu j}^n(\xi)=\mu_{\nu j}(\xi)+L_{\nu}(\xi) \mu_{0,j}(\xi),\\
	\w\mu_{0,j}^n(\xi)&=\w\mu_{0,j}(\xi)-
	\sum\limits_{i=1}^r\overline{L_i}(\xi)\w\mu_{i,j}(\xi), \quad
	& &\w\mu_{\nu j}^n(\xi)=\w\mu_{\nu j}(\xi). 
	\end{alignedat}
	\end{equation*}       
	Let new masks $m_{\nu}^n$, $\w m_{\nu}^n$ be constructed from the new  polyphase components by~(\ref{PR}). Thus, the transformation formulas for the masks  are
	\begin{equation}
	\begin{alignedat}{2}
	m_0^n(\xi)&=m_0(\xi), \quad
	& \quad & m_{\nu}^n(\xi)=m_{\nu}(\xi)+L_{\nu}(M^*\xi) m_0(\xi),
	\\
	\w m_0^n(\xi)&= \w m_0(\xi)-\sum\limits_{i=1}^r\overline{L_i}(M^*\xi)\w m_i(\xi), \quad 
	& &\w m_{\nu}^n(\xi) =\w m_{\nu}(\xi).
	\label{fLiftTransform}
	\end{alignedat}
	\end{equation} 
	The above transformation of the masks is called  the \textit{lifting
		scheme transformation}.
	Note that the choice of trigonometric polynomials 
	$L_{\nu},$ ${\nu}=1,\dots,r$ 
	is not restricted.

The aim is to find the  lifting
scheme transformation such that all new masks 
preserve their symmetry properties	and	
all new wavelet masks have vanishing moments at least of order 1.

\begin{theo}
	Let a dilation matrix $M$ and a center
	$c$ be  appropriate for a symmetry group ${\cal H}$, $n\in\n$.
	Suppose $m_0$ and $\w m_0$ are as in Theorem~\ref{theoWaveNoSym}.
	Suppose $m_{(p,i)},$ $\w m_{(p,i)}$ are wavelet masks constructed using
	Theorem~\ref{theoWaveNoSym}.
	Assume that trigonometric polynomials $L_{p,i},$  satisfy
	$
	L_{p,i}(\nul)=-m_{(p,i)}(\nul)
	$ and
	\be
	L_{p,0}(M^*F^*M^{*-1}\xi)= L_{p,0}(\xi)e^{2\pi i (r_{p,0}^F,\xi)},
	\quad \forall F \in {\cal H}_{p,0},
	\label{fLFsym}
	\ee
	\be
	L_{p,0}(M^*E^{(i)*}M^{*-1}\xi)= L_{p,i}(\xi),
	\quad E^{(i)} \in {\cal E}_p,   
	\label{fLEsym}
	\ee
	for $i=0,\dots,\#{\cal E}_{p}-1,$ $p=0,\dots,\#\Lambda-1$. New masks $m_0^n,$ $\w m_0^n$, $m_{(p,i)}^n$, $\w m_{(p,i)}^n$ defined by the lifting
	scheme transformation preserve symmetry properties of masks
	$m_0,$ $\w m_0$, $m_{(p,i)}$, $\w m_{(p,i)}$, respectively.
	New wavelet masks $m_{(p,i)}^n$ have vanishing moments at least of order 1.
	If new refinable functions $\phi,$ $\w\phi$
	corresponding to new refinable masks $m^n_0$, $\w m^n_0$
	are in $L_2(\rd)$, then the resulting wavelet system
	is a dual wavelet frame.
	\label{theoLiftSym}
\end{theo}

{\bf Proof.}		
Vanishing moments at least of order 1 for wavelet masks
$m_{(p,i)}^n$, are provided
by conditions
$L_{p,i}(\nul)=-m_{(p,i)}(\nul).$
This follows from~(\ref{fLiftTransform}).

Next, we show that the symmetry properties are preserved.
For $m_{(p,0)}^n$ condition~\ref{fLS_1cond} in Theorem~\ref{theoWaveNoSym} is preserved, since
conditions~(\ref{fLiftTransform}) and~(\ref{fLFsym}) yield
\begin{equation}
\begin{split}
m_{(p,0)}^n(F^*\xi) &=
m_{(p,0)}(F^*\xi)+L_{(p,0)}(M^*F^*\xi) m_0(F^*\xi) \\
& = m_{(p,0)}(\xi)e^{2\pi i (Fs_{p,0}-s_{p,0},\xi)}+
L_{(p,0)}(M^*F^*\xi) m_0(\xi)e^{2\pi i (Fc-c,\xi)} \\
& =  m_{(p,0)}^n(\xi)e^{2\pi i (Fs_{p,0}-s_{p,0},\xi)}.
\label{eq:LiftSymPreserve}
\end{split}
\end{equation} 
Next, condition~\ref{fLS_2cond}  in Theorem~\ref{theoWaveNoSym} for $m^n_{(p,i)}$ is preserved, since
conditions~(\ref{fLiftTransform}) and~(\ref{fLEsym}) yield  
\begin{equation}
\begin{split}
m_{(p,0)}^n(E^{(i)*}\xi)& =
m_{(p,0)}(E^{(i)*}\xi)+L_{(p,0)}(M^*E^{(i)*}\xi) m_0(E^{(i)*}\xi) \\
& = m_{(p,i)}(\xi)e^{2\pi i (E^{(i)}c-c,\xi)}+
L_{(p,i)}(M^*\xi)
m_0(\xi)e^{2\pi i (E^{(i)}c-c,\xi)} \\ & =
m_{(p,i)}^n(\xi)e^{2\pi i (E^{(i)}c-c,\xi)},
\label{eq:LiftSymPreserve2}
\end{split}
\end{equation} 
where 
$E^{(i)}\in{\cal E}_{p},$
$i=0,\dots,\#{\cal E}_{p}-1,
$
$p=0,\dots,\#\Lambda-1.$    

Now, we show that new refinable mask $\w m^n_0$ remains ${\cal H}$-symmetric with respect to the point $c$. 
Let us fix $p.$
Recall that ${\cal H}$ can be uniquely represented as follows
${\cal H} = {\cal E}_p \times {\cal H}_{p,0}.$
Suppose $K \in {\cal H}$ and $E^{(i)} \in {\cal E}_p$. Define a mapping $j(\cdot,p,K)$ from the set of indices $ \{0,\dots ,\#{\cal E}_p- 1\}$ to
itself, where index $j = j(i,p,K)$ is uniquely defined such that $KE^{(i)}=E^{(j)} F$. 

For $K\in {\cal H}$, consider
\be
\w m^n_0(K^*\xi)=\w m_0(\xi) e^{2\pi i (Kc-c,\xi)}-
\sum\limits_{p=0}^{\#\Lambda-1} 
\sum\limits_{i=0}^{\#{\cal E}_{p}-1}
\overline{L_{p,i}(M^*K^*\xi)} \w m_{p,i}(K^*\xi)  .
\label{fLSCheck}
\ee

Let us fix indices $p$ and $i.$ Then
$KE^{(i)}=E^{(j)} F$, where $F\in {\cal H}_{p,0}$,
$E^{(j)}\in {\cal E}_{p}$, $j=j(i,p,K)$.
With conditions on $L_{p,i},$ we obtain
\begin{equation*}
\begin{split}
\overline{L_{p,i}(M^*K^*\xi)}=\overline{L_{p,0}(M^*E^{(i)*}K^*\xi)}& =
\overline{L_{p,0}(M^*F^*E^{(j)*}\xi)}\\
&=\overline{L_{p,0}(M^*E^{(j)*}\xi)} 
e^{-2\pi i (M r_{p,0}^F,E^{(j)*}\xi)}\\ &=
\overline{L_{p,j}(M^*\xi)}e^{-2\pi i (E^{(j)}M r_{p,0}^F,\xi)},
\end{split}
\end{equation*} 
where 
$-E^{(j)}M r_{p,0}^F= E^{(j)}s_{p,0}-E^{(j)}c-KE^{(i)}(s_{p,0}- c).$
Using~\ref{fLS_1cond} and~\ref{fLS_2cond} we get
\begin{equation*}
\begin{split}
\w m_{(p,i)}(K^*\xi)&=
\w m_{(p,0)}(E^{(i)*}K^*\xi)
e^{2\pi i (c-E^{(i)}c,K^*\xi)}\\
&=
\w m_{(p,0)}(F^*E^{(j)*}\xi)  e^{2\pi i (c-E^{(i)}c,K^*\xi)}\\ 
& =	\w m_{(p,0)}(E^{(j)*}\xi)
e^{2\pi i (Fs_{p,0}-s_{p,0},E^{(j)*}\xi)}
e^{2\pi i (c-E^{(i)}c,K^*\xi)} \\
& =	 \w m_{(p,j)}(\xi)
e^{2\pi i (E^{(j)}c-c,\xi)}
e^{2\pi i (E^{(j)}(Fs_{p,0}-s_{p,0}),\xi)}
e^{2\pi i (Kc-KE^{(i)}c,\xi)}.
\end{split}
\end{equation*} 
So, finally we have
$$
\overline{L_{p,i}(M^*K^*\xi)}\w m_{(p,i)}(K^*\xi)=
\w m_{(p,j)}(\xi)\overline{L_{p,j}(M^*\xi)}e^{2\pi i (R,\xi)},
$$
where 
\begin{equation*}
\begin{split}
R&=-E^{(j)}M r_{p,0}^F+E^{(j)}c-c+E^{(j)}(Fs_{p,0}-s_{p,0})+Kc-KE^{(i)}c\\ & =
E^{(j)}s_{p,0}-E^{(j)}c-KE^{(i)}(s_{p,0}- c)+
E^{(j)}c-c +
\\ &~~~~~~~~~~~~~~~~~~~~~~~~~~~~~~~~~~~~~~~~~~KE^{(i)}s_{p,0}-E^{(j)}s_{p,0}+Kc-KE^{(i)}c\\ & = Kc-c.
\end{split}
\end{equation*} 
Thus, by~(\ref{fLSCheck}) we obtain that $\w m_0^n$ is ${\cal H}$-symmetric
with respect to the center $c$. $\Diamond$

Condition~(\ref{fLFsym}) means that $L_{p,0}$ is $M^{-1} H_{p,0} M$--symmetric with respect to the center $M^{-1}(s_{p,0}-c)$, $p=0,\dots, \#\Lambda-1.$ 
Such trigonometric polynomials satisfying condition $L_{p,0} (\nul)= - m_{(p,0)}(\nul)$ can be easily constructed. 

\subsection{Symmetric wavelet frames}

In~\cite[Algorithm 1]{Sk1} an algorithm for the construction of
dual wavelet frames was suggested. Here we modify this algorithm,
such that the constructed wavelets have the symmetry properties.

Firstly, we give  a constructive description of the method. 
Let 
a dilation matrix $M$ and a center $c$ be appropriate for ${\cal H}$,
$n\in\n$.

\centerline{{\bf Algorithm 1}}

\noindent{\bf Step 1.} 
Find  a refinable mask $m_0$ that is ${\cal H}$-symmetric with respect 
to the center $c$ and obeys the sum rule of order $n$. Let $(2\pi i)^{|\beta|} \lambda_{\beta}:=D^\beta  m_0(\nul)$, $\beta\in\Delta_n.$

\noindent{\bf Step 2.} 
Next, define numbers $\w\lambda_\beta$, $\beta\in\Delta_n,$ such that~(\ref{fLambda}) is valid and 
find a trigonometric polynomial $\w m'_0$ that is ${\cal H}$-symmetric with respect 
to the center $c$, obeys the sum rule of order $n$ and $(2\pi i)^{|\beta|} \w \lambda_{\beta}:=D^\beta \w m'_0(\nul)$.

\noindent{\bf Step 3.}
Let $\widetilde\mu'_{00}\ddd \widetilde\mu'_{0,m-1}$ be the polyphase components of $\w m'_0$.
Set 
\ba
\sigma:=\sum_{l=0}^{m-1}\overline{\mu_{0l}}\,\widetilde\mu'_{0l},\quad
\w\mu_{0k}&:=&{(2-\sigma)}\w\mu^\prime_{0k}, k=0\ddd m-1,
\label{fSigma}
\\
\mu_{0m}&:=&\overline{(1-\sigma)},\  \widetilde\mu_{0m}:={1-\sigma}.
\label{fNewMu}
\ea
Due to~\cite[Corollary 15]{Sk1}, we have $D^\beta\sigma(\nul)=0$
for all   $\beta\in\zd_+$,  $0<|\beta|< n$. Therefore, 
$D^\beta\lll \overline{1-\sigma}\rrr(\nul)=0$
for all   $\beta\in\zd_+$,  $0\le|\beta|< n$.  It follows that
$D^\beta\w\mu_{0k}(\nul)=D^\beta\w\mu^\prime_{0k}(\nul)$
for all  $\beta\in\zd_+$,  $|\beta|< n$.
Therefore,
$$
1-\sum_{k=0}^{m-1}\overline{\mu_{0k}}\widetilde\mu_{0k}=1-(2-\sigma)
\sum_{k=0}^{m-1}
\overline{\mu_{0k}}\w\mu'_{0k}
=(1-\sigma)^2,
$$
which yields $\sum_{k=0}^{m}\overline{\mu_{0k}}\widetilde\mu_{0k}=1$.
Dual refinable mask $\w m_0$ is constructed by~(\ref{PR}) from $\w\mu_{0k}$.
Note that $\w m_0$ obeys the sum rule of order $n$. 
Since $\sigma$ is ${\cal H}$-symmetric with respect to the origin (the proof is below in Theorem~\ref{theoFrameNoSym}), then dual refinable mask   $\w m_0$ is ${\cal H}$-symmetric with respect to the center $c$.

\noindent{\bf Step 4.}
Set  $r=m$ if $\sigma\equiv1$, otherwise we set $r=m+1$, $\mu_{0,m+1}\equiv0, \widetilde\mu_{0, m+1}\equiv0$. Let ${\cal N}:=\{\mu_{\nu k}\}_{\nu,k=0}^{r}$,
$\widetilde{\cal N}:=\{{\widetilde \mu_{\nu k}}\}_{\nu,k=0}^{r}$.
Then matrix extension  
can be realized as follows
\be
{\cal N}:=\left(%
\begin{array}{cc}
	P & 0 \\
	I_{m} - \w P^* P & \w P^*\\
\end{array}%
\right),\ \ \
\widetilde{\cal N}:=\left(%
\begin{array}{cc} 
	\w P & 0 \\
	I_{m} - P^* \w P & P^*\\
\end{array}%
\right), 
\label{eq:MEPframe1}
\ee
if $r=m$, or
\be
{\cal N}:=\left(%
\begin{array}{ccc}
	P & \mu_{0m} & 0 \\
	I_{m} - \w P^* P & \mu_{0m} \w P^* & \w P^*\\
	- \overline{\w\mu_{0m}} P &  1 - \overline{\w\mu_{0m}} \mu_{0m} & \overline{\w\mu_{0m}}\\
\end{array}%
\right),\ \ \
\widetilde{\cal N}:=\left(%
\begin{array}{ccc} 
	\w P & \w\mu_{0m} & 0 \\
	I_{m} - P^* \w P & \w \mu_{0m} P^* & P^*\\
	- \overline{\mu_{0m}} \w P &  1 - \overline{\mu_{0m}} \w\mu_{0m} & \overline{\mu_{0m}}
\end{array}%
\right), 
\label{eq:MEPframe2}
\ee
if $r=m+1$,
where $P=(\mu_{00},\dots,\mu_{0,m-1})$,
$\w P=(\w \mu_{00}, \dots,\w \mu_{0,m-1}).$
It is not difficult to see
that the  matrices 
satisfy ${\cal N}{\widetilde {\cal N}^*}=I_{r}$. This yields that the columns of the polyphase matrices $\cal M, \widetilde {\cal M}$ are biorthonormal.
Wavelet masks 
$m_\nu, \w m_\nu$, $\nu=1\ddd r$, are constructed from the polyphase components.~ $\Diamond$

The next Theorem states that wavelet masks constructed by Algorithm 1 have
the mutual symmetry properties.
Again, for convenience, we enumerate wavelet masks using double index 
as in Theorem~\ref{theoWaveNoSym}.
For the submatrices $I_{m} - \w P^* P$ and $I_{m} - P^* \w P$
we use notations as in~(\ref{eq:Nblock}).

\begin{theo} 
	Let a dilation matrix $M$ and a center
	$c$ be  appropriate for a symmetry group ${\cal H}$, $n\in\n$.
	Suppose refinable and wavelet masks 
	are constructed by
	Algorithm 1.
	Then wavelet masks  $m_{(p,i)}$ and $\w m_{(p,i)}$
	have the following symmetry properties:
	\begin{itemize}
		\item $m_{(p,0)}$, $\w m_{(p,0)}$ are ${\cal H}_{p,0}$-symmetric 
		with respect to   the center $s_{p,0};$ 
		\item $m_{(p,i)}=m_{(p,0)}(E^{(i)*}\xi)e^{2\pi i (c-E^{(i)}c,\xi)},$
		$\w m_{(p,i)}=\w m_{(p,0)}(E^{(i)*}\xi) e^{2\pi i (c-E^{(i)}c,\xi)},$  
		$E^{(i)}\in{\cal E}_{p},$
	\end{itemize}    
	$i=0,\dots,\#{\cal E}_{p}-1,$
	$p=0,\dots,\#\Lambda-1$.
	If $r=m+1$, then wavelet masks $m_r$ and $\w m_r$ defined by 
	the last rows of matrices ${\cal N}$ and $\w {\cal N}$ are 
	${\cal H}$-symmetric with respect to the center $c$.
	All wavelet masks have vanishing moments of order  $n$.
	If refinable functions $\phi,$ $\w\phi$
	corresponding to refinable masks $m_0$, $\w m_0$
	are in $L_2(\rd)$, then the resulting wavelet system
	is a dual wavelet frame.
	\label{theoFrameNoSym}
\end{theo}

{\bf Proof.}
First, we prove that $\sigma$ defined in~(\ref{fSigma}) is ${\cal H}$-symmetric with respect to the origin.
Let us rewrite $\sigma$ using another enumeration of the polyphase components 
and show that $\sigma(M^* K^* M^{*-1}\xi) = \sigma(\xi)$ for all $K \in {\cal H}.$ Thus,
\begin{equation*}
\begin{split}
\sigma(M^* K^* M^{*-1}\xi)  = & \sum_{p=0}^{\#\Lambda-1}
\sum_{i=0}^{\#{\cal E}_p-1} \overline{\mu_{0,p,i}(M^* K^* M^{*-1})}\,\widetilde\mu'_{0,p,i}(M^* K^* M^{*-1})  \\  = &
\sum_{p=0}^{\#\Lambda-1}
\sum_{i=0}^{\#{\cal E}_p-1} \overline{\mu_{0,p,0}(M^* E^{(i)*} K^* M^{*-1}\xi)}\,\widetilde\mu'_{0,p,0}(M^* E^{(i)*} K^* M^{*-1}\xi).
\end{split}
\end{equation*}
Let us fix $p.$
Recall that ${\cal H}$ can be uniquely represented as follows
${\cal H} = {\cal E}_p \times {\cal H}_{p,0}.$
Recall that for a matrix $K$ and fixed $p$ 
there exist unique matrices  
$E^{(j)}\in {\cal E}_{p}$  and 
$F\in {\cal H}_{p,0}$ such that 
$KE^{(i)}=E^{(j)} F$, where $j = j(p,i,K)$.
Using~(\ref{fPolyF}),~(\ref{fPolyE_i}) and continuing the above equalities,
we get
\begin{equation*}
\begin{split}
\sigma(M^* K^* M^{*-1}\xi) = & \sum_{p=0}^{\#\Lambda-1}
\sum_{i=0}^{\#{\cal E}_p-1} \overline{\mu_{0,p,0}(M^* F^* E^{(j(p,i,K))*}  M^{*-1}\xi)}\,\widetilde\mu'_{0,p,0}(M^*  F^* E^{(j(p,i,K))*}  M^{*-1}\xi)
\\
= & \sum_{p=0}^{\#\Lambda-1}
\sum_{i=0}^{\#{\cal E}_p-1} \overline{\mu_{0,p,j(p,i,K)}(\xi)}\,\widetilde\mu'_{0,p,j(p,i,K)}(\xi) = \sigma(\xi).
\end{split}
\end{equation*}
Since $\sigma$ is ${\cal H}$-symmetric with respect to the origin,
then polyphase components $\w\mu'_{0k}$ and $\w\mu_{0k}$ defined in~(\ref{fNewMu}) have the same symmetry properties. Therefore, dual refinable mask $\w m_0$ is ${\cal H}$-symmetric with respect to the center $c.$
Symmetry properties of wavelet masks $m_{(p,i)}$, $\w m_{(p,i)}$
defined by submatrices $I_{m} - \w P^* P$ and $I_{m} - P^* \w P$
in~(\ref{eq:MEPframe1}) or~(\ref{eq:MEPframe2})
can be proved analogously to the proof of Theorem~\ref{theoWaveNoSym}.
If $r=m+1$, then the last wavelet masks $m_r$ and $\w m_r$ defined by 
the last rows of matrices ${\cal N}$ and $\w {\cal N}$ are 
${\cal H}$-symmetric with respect to the center $c$.
Indeed, since
$$
m_{r} (\xi) = - \overline{\w \mu_{0m}(M^* \xi)} m_0(\xi),
\quad
\w m_{r} (\xi) = - \overline{\mu_{0m}(M^* \xi)} \w m_0(\xi),
$$
then for all $K\in {\cal H}$
\begin{equation*}
\begin{split}
m_{r} (K^*\xi) & = - \overline{\w \mu_{0m}(M^*K^*M^{*-1} M^* \xi)} m_0(K^*\xi) \\ & =  - \overline{\w \mu_{0m}(M^* \xi)} m_0(\xi) e^{2\pi i (Kc-c,\xi)} \\ & = m_{r} (\xi) e^{2\pi i (Kc-c,\xi)}
\end{split}
\end{equation*}
and analogously with $\w m_{r}$.

Vanishing moments for all wavelet masks are provided by~\cite[Theorem 8]{Sk1}.
$\Diamond$

Note that wavelet masks $m_{(p,i)}$, $\w m_{(p,i)}$ are ${\cal H}_{p,i}$--symmetric  with respect to   the center $s_{p,i}$. 


The resulting wavelet masks constructed by Theorem~\ref{theoFrameNoSym} can have quite huge coefficient supports. It can be reduced if we will not require vanishing moments of order $n$ for all wavelet functions. According to Theorem 7 in~\cite{Sk1}, the vanishing moments of order $n$ for dual wavelet masks are provided by sum rule of order $n$ for $m_0$ and  $D^\beta\lll{1-\sigma}\rrr(\nul)=0$
for all  $|\beta|< n$. These vanishing moments should be kept, since they provide approximation order $n$ for the resulting wavelet system. The vanishing moments of order $n$ for wavelet masks are provided by sum rule of order $n$ for $\w m_0$ and $D^\beta\lll{1-\sigma}\rrr(\nul)=0$
for all  $|\beta|< n$. Here, the order of vanishing moments can be reduced up to $1$ (since we need only to keep the frame condition). So, we can require the sum rule of order $1$ for $\w m_0'$, preserving $D^\beta\lll {1-\sigma}\rrr(\nul)=0$ for all  $|\beta|< n$. In some cases, this can be done (see Example 2 in Section 4).

\section{Symmetrization}

The aim of this section is to construct  wavelet masks such that all of them are ${\cal H}$-symmetric in a sense. 
To
do that we extend the definition of ${\cal H}$-symmetric trigonometric polynomials.
Let $t$ be a trigonometric polynomial. Then $t$ has the ${\cal H}$-symmetry property if for each matrix $E\in {\cal H}$ 
$$t(E^*\xi)= \varepsilon_E
e^{2\pi i (r_E,\xi)}         
t(\xi),$$
where $\varepsilon_E\in {\mathbb C}$, $|\varepsilon_E|=1,$ 
$r_E\in {\mathbb Z}^{d}.$ 
This definition is a generalization of~(\ref{TrigSym}).
For example, the new definition includes antisymmetric trigonometric polynomials and trigonometric polynomials, which do not have a symmetry center common to all matrices $E\in {\cal H}$.

Let ${\cal H}$ be an abelian symmetry group.                
Assume we have a row of the polyphase components of an ${\cal H}$-symmetric mask.  Next, we find a unitary transformation of the row such that each element of the new row has the ${\cal H}$-symmetry property. 
For a fixed $p$, ${\cal H}_{p,0}$ and ${\cal E}_{p}$ are an abelian subgroups of ${\cal H}$ and ${\cal E}_{p}$ can be expressed as the direct product of cyclic subgroups. Let $\gamma_p$ be the number of cyclic subgroups of ${\cal E}_{p}$.
Then there exist unique prime numbers $N_{p,j},$ $j=1,\dots,\gamma_p$ and matrices ${\cal K}_1,\dots,{\cal K}_{\gamma_p}\in {\cal E}_{p}$ such that
$${\cal E}_{p}=\{I_d,{\cal K}_1,\dots,{\cal K}_1^{N_{p,1}-1}\}\times\dots
\times\{I_d,{\cal K}_{\gamma_p},\dots,{\cal K}_{\gamma_p}^{N_{p,\gamma_p}-1}\}.$$

Thus, fixed element $E\in {\cal E}_{p}$ can be uniquely
represented as follows $E=\prod_{j=1}^{\gamma_p} {\cal K}_j^{k_j}$, where $k_j$ are some integers from the sets 
$\{0,\dots,N_{p,j}-1\},$ respectively.
So, this matrix $E$ we denote by $E^{(k)}$ where $k=(k_1,\dots,k_{\gamma_p})$ is a number in a mixed radix number system for set of numbers
$\{0,\dots,\#{\cal E}_{p}-1\}$ with the base $(N_{p,1},\dots,N_{p,\gamma_p})$.   
For two matrices $E^{(k)}, E^{(l)}\in {\cal E}_{p}$ their product 
can be written as follows  $E^{(k)}E^{(l)}=E^{(k \oplus l)}$, where 
$k \oplus l$ is an addition on the set  $\{0,\dots,\#{\cal E}_{p}-1\}$ defined as follows:
$$k\oplus l=((k_1+l_1)\mod N_{p,1},\dots,(k_{\gamma_p}+l_{\gamma_p})\mod N_{p,\gamma_p}).$$

For the integers $N_{p,j}$ we denote
$\varepsilon_{N_{p,j}}:= e^{ \frac {2\pi i}{N_{p,j}}}.$    
For any $p\in\{1,\dots,\#\Lambda-1\}$,
define
the matrix of the discrete Fourier transform
$W_{N_{p,j}}=\frac 1 {\sqrt{N_{p,j}}}\{
\varepsilon_{N_{p,j}}^{kl}\}_{k,l=0,N_{p,j}-1}$. It is known that $W_{N_{p,j}}$ is a unitary and symmetric matrix, i.e. $W_{N_{p,j}}W_{N_{p,j}}^*=I_m, W_{N_{p,j}}^T=W_{N_{p,j}}.$
Define ${\cal W}_p=W_{N_{p,1}}\otimes\dots\otimes W_{N_{p,\gamma_p}},$ where operation $\otimes$ is the Kronecker product.
${\cal W}_p$ is a $(\#{\cal E}_{p})\times(\#{\cal E}_{p})$ unitary matrix. Some properties of the matrix ${\cal W}_p$ are $$\left[{\cal W}_p\right]_{k,l}
\left[{\cal W}_p\right]_{n,l}=\left[{\cal W}_p\right]_{k\oplus n,l}, \quad \left[{\cal W}_p\right]_{k,l}\overline{\left[{\cal W}_p\right]_{k,l}}=1, \quad 
k,l,n=0,\dots,\#{\cal E}_{p}-1.$$

The next Lemma states that ${\cal W}_p$ symmetrizes the part of the row of polyphase components.

\begin{lem}~\cite[Lemma 17]{KrSym} Suppose $m_0$ is an ${\cal H}$-symmetric with respect to the an appropriate center $c$ mask and $\mu_{0,p,i}$ are its polyphase components. For a fixed $p$, a row $T_p$ is  defined by 
	$T_p=(\mu_{0,p,0},...,\mu_{0,p,\#{\cal E}_{p}-1}).$
	Suppose that $r_{p,0}^F=M^{-1}E Mr_{p,0}^F$ for all  $F\in {\cal H}_{p,0}$ and $E\in{\cal E}_{p}.$ Then each element of the row $T_p':=T_p{\cal W}_p$ has the ${\cal H}$-symmetry property, i.e. for any $K\in {\cal H}$ 
	\begin{equation}
	\mu'_{0,p,r}((M^{-1}K M)^*\xi)=
	\overline{\left[{\cal W}_p\right]_{i,r}}\mu'_{0,p,r}(\xi)e^{-2\pi i (r_{p,0}^{F},\xi)},
	\quad r=0,\dots,\#{\cal E}_p-1,
	\label{fHSymPropRow}
	\end{equation}
	where $K=E^{(i)} F$, $F\in {\cal H}_{p,0}$ and $E^{(i)}\in{\cal E}_{p},$
	$T_p' = (\mu'_{0,p,0},...,\mu'_{0,p,\#{\cal E}_{p}-1})$
	\label{propSym}
\end{lem}

According with this Lemma we need a \textit{special assumption  } \be
r_{p,0}^F=M^{-1}E Mr_{p,0}^F, \quad \forall F\in {\cal H}_{p,0}, \forall E\in{\cal E}_{p}
\label{fSpecAssum}
\ee   to ensure the ${\cal H}$-symmetry property for all components of the row $T_p'$.
Notice that this assumption is used only to provide $M^{-1}{\cal H}_{p,0}M$-symmetry for $\mu'_{0,p,r}(\xi)$.

\subsection{Symmetrization for frame-like wavelets}

Define a block diagonal unitary matrix ${\cal W}$ as follows:
${\cal W}=\texttt{diag} ({\cal W}_{0},\dots,{\cal W}_{\#\Lambda-1}).$ Let $m_0$ be an ${\cal H}$-symmetric mask.
Suppose~(\ref{fSpecAssum}) is valid for all 
$p=0,\dots,\#\Lambda-1.$
Then by Lemma~\ref{propSym} matrix ${\cal W}$ symmetrizes the row $P=(T_0,\dots,T_{
	\#\Lambda-1})$ of the polyphase components of $m_0,$ namely all elements of the row 
$$P^s:=P{\cal W}=(T_0{\cal W}_0,\dots,T_{\#\Lambda-1}{\cal W}_{\#\Lambda-1})$$
have the ${\cal H}$-symmetry property. 
Wavelets with the ${\cal H}$-symmetry property can be constructed using matrix ${\cal W}$.   

\begin{theo} 
	Let a dilation matrix $M$ and a center
	$c$ be  appropriate for an abelian symmetry group ${\cal H}$, $n\in\n$.
	Suppose $m_0$ and $\w m_0$ are ${\cal H}$-symmetric
	with respect to the center $c$ masks such that
	mask $m_0$
	obeys the sum rule of order $n$,
	mask $\w m_0$ satisfies
	condition~(\ref{20_new}).
	Suppose that condition~(\ref{fSpecAssum}) is valid for all 
	$p=0,\dots,\#\Lambda-1.$
	Then there exist wavelet masks 
	$m_{(p,i)}$ and $\w m_{(p,i)}$  
	which have the ${\cal H}$-symmetry property,
	wavelet masks $\w m_{(p,i)}$ have vanishing moments of order $n$,
	$i=0,\dots,\#{\cal E}_{p}-1,$
	$p=0,\dots,\#\Lambda-1$.
	The corresponding MRA-based dual wavelet system 
	is almost frame-like in $S'.$   
	\label{theoWaveSym}
\end{theo}

{\bf Proof.}   
Let $\mu_{0,p,i}$ and $\w\mu_{0,p,i}$ be the polyphase components of $m_0$ and $\w m_0$, 
$i=0,\dots,\#{\cal E}_{p}-1,$     $p=1,\dots,\#\Lambda-1,$
and let 
$T_p=(\mu_{0,p,0},...,\mu_{0,p,\#{\cal E}_p-1}),$ and 
$\w T_p=(\w\mu_{0,p,0},...,\w\mu_{0,p,\#{\cal E}_p-1}).$
Set $P=(T_0,\dots,T_{\#\Lambda-1}),$ 
$\w P=(\w T_0,\dots,\w T_{\#\Lambda-1})$.
Let us consider matrix extension~(\ref{fTTblock}) with $U=\w U={\cal W}^*$
\be
{\cal N}=\left( \begin{array}{c|c}
	P & 1-P\w P^* \\
	\hline
	&  \\
	{\cal W}^* & -{\cal W}^*\w P^{*}
\end{array}\right),
\quad 
{\cal \w N}=\left(
\begin{array}{c|c}
	\w P& 1 \\
	\hline
	&   \\
	{\cal W}^*-{\cal W}^* P^{*} \w P & -{\cal W}^*P^{*}
\end{array} \right).
\label{fNNblockU}
\ee

Let us consider submatrix $ {\cal W}^*-{\cal W}^* P^{*} \w P$.
It is a block matrix
$$ \begin{footnotesize}
{\cal W}^*-{\cal W}^* P^{*} \w P=
\left( \begin{array}{cccc}
{\cal W}_0^*-(T_0 {\cal W}_0)^*\w T_0 & 
-(T_0 {\cal W}_0)^*\w T_1 & \dots & 
-(T_0 {\cal W}_0)^*\w T_{\#\Lambda-1}\\
-(T_1 {\cal W}_1)^*\w T_0 &
{\cal W}_1^*-(T_1 {\cal W}_1)^*\w T_1 & \dots & 
-(T_1 {\cal W}_1)^*\w T_{\#\Lambda-1} \\
\vdots & & \ddots & \vdots \\
-(T_{\#\Lambda-1} {\cal W}_{\#\Lambda-1})^*\w T_0 &	\dots &\dots &
{\cal W}_{\#\Lambda-1}^*-(T_{\#\Lambda-1} 
{\cal W}_{\#\Lambda-1})^*\w T_{\#\Lambda-1}
\end{array}\right).
\end{footnotesize}
$$
Denote the elements of the submatrix as 
$ {\cal W}^*-{\cal W}^* P^{*} \w P=\{\w\mu_{(p,i),(t,j)}\}_
{p=0,\dots,\#\Lambda-1, i=0,\dots,\#{\cal E}_p-1}^
{t=0,\dots,\#\Lambda-1, j=0,\dots,\#{\cal E}_t-1}.$
With this enumeration the element $\w\mu_{(p,i),(t,j)}$ in the submatrix is the element in the block $(p,t)$ with position $(i,j)$. Namely, if $p\neq t$ then
$$\w\mu_{(p,i),(t,j)}(\xi)=[-(T_p {\cal W}_p)^*\w T_t]_{i,j}=
-\overline{\mu'_{0,p,i}(\xi)}\,\,\w\mu_{0,t,j}(\xi);$$
if $p=t$ then
$$\w\mu_{(p,i),(p,j)}(\xi)=[ {\cal W}^*_p-(T_p {\cal W}_p)^*\w T_p]_{i,j}=
\overline{[ {\cal W}_p]_{i,j}}-\overline{\mu'_{0,p,i}(\xi)}\,\,\w\mu_{0,p,j}(\xi),$$
where $\mu'_{0,p,i},$ $ i=0,\dots,\#{\cal E}_p-1$, are the elements of the row $T_p {\cal W}_p$,
$ p=0,\dots,\#\Lambda-1$. All $\mu'_{0,p,i}$ have the ${\cal H}$-symmetry property by Lemma~\ref{propSym}.

Next, we collect 
wavelet masks by~(\ref{PR}). For fixed  $ p=0,\dots,\#\Lambda-1$,
$i=0,\dots,\#{\cal E}_{p}-1,$ we get
$$
\w m_{(p,i)}(\xi)=\frac{1}{\sqrt{m}}\sum_{j=0}^{\#{\cal E}_{p}-1} \overline{[W_p]_{i,j}}
e^{2\pi i (s_{p,j},\xi)}-\overline{\mu'_{0,p,i}(M^{*}\xi)}\w  m_0(\xi).
$$
Check that $\w m_{(p,i)}$ has the ${\cal H}$-symmetry property. 
For a fixed $F\in {\cal H}_{p,0}$, by~(\ref{fHSymPropRow}), 
and by the properties of ${\cal W}_p$ we obtain
\begin{equation*} 
	\begin{split}
\w m_{(p,i)}(F^*\xi)=&\frac{1}{\sqrt{m}}\sum_{j=0}^{\#{\cal E}_{p}-1} \overline{[W_p]_{i,j}}
e^{2\pi i (F s_{p,j},\xi)}-\overline{\mu'_{0,p,i}(M^{*}F^*\xi)}\w m_0(F^*\xi) \\
=   &  
e^{2\pi i (F s_{p,0}-s_{p,0},\xi)}\w m_{(p,i)}(\xi),
	\end{split}
\end{equation*} 
since $F s_{p,j} = F (E^{(j)} s_{p,0} + c - E^{(j)} c) = 
s_{p,j} + F s_{p,0} - s_{p,0}$ by~(\ref{fSpecAssum}),~(\ref{fDigF}),
(\ref{fDigE_i}).
For a fixed $E^{(k)}\in {\cal E}_p,$ by~(\ref{fHSymPropRow}), 
and by the properties of ${\cal W}_p$ we obtain
\begin{equation*}
\begin{split} 
\w m_{(p,i)}(E^{(k)}\xi)=& \frac{1}{\sqrt{m}}\sum_{j=0}^{\#{\cal E}_{p}-1} \overline{[W_p]_{i,j}}
e^{2\pi i (E^{(k)} s_{p,j},\xi)}-\overline{\mu'_{0,p,i}(M^{*}E^{(k)}\xi)}\w m_0(E^{(k)}\xi) 
\\   
= & [{\cal W}_p]_{k,i} e^{2\pi i (E^{(k)}c - c,\xi)}\w m_{(p,i)}(\xi),
\end{split}
\end{equation*}
since $E^{(k)} s_{p,j} = 
s_{p,j\oplus k} + E^{(k)}c - c$ by~(\ref{fDigE_i}). 
Therefore, for a fixed $K\in {\cal H}$, with 
$K=E^{(k)}F$, $E^{(k)}\in {\cal E}_{p}$,
$F\in {\cal H}_{p,0},$ 
$$
\w m_{(p,i)}(K^*\xi) = 
[W_p]_{k,i}  e^{2\pi i (Kc - c + M r_{p,0}^F,\xi)}\w m_{(p,i)}(\xi).
$$

Thus, we get the ${\cal H}$-symmetry property for 
the wavelet mask $\w m_{(p,i)}$.
It is not hard to see that the wavelet masks $ m_{(p,i)}$ have the same ${\cal H}$-symmetry property as $\w m_{(p,i)}$.

The vanishing moments of order $n$ for the  wavelet masks $\w m_{(p,i)}$ are provided by Theorem~\ref{lemKrSk}.              
$\Diamond$	

Note that if~(\ref{fSpecAssum}) is not valid for some $p$, then we can take the corresponding matrix ${\cal W}_p$ equal to the identity matrix. Thus, for this index $p$ the corresponding set of wavelet masks  $m_{(p,i)}$, $\w m_{(p,i)}$ will remain mutually symmetric.

\subsection{Symmetrization for lifting scheme}

In this subsection we suggest a lifting scheme transformation preserving
the ${\cal H}$-symmetry property of the constructed by Theorem~\ref{theoWaveSym} wavelet masks.
The corresponding trigonometric polynomials for the transformation we defined as follows: let ${\cal L}_p = (L_{p,0},\dots,L_{p,\#{\cal E}_p-1}).$
Define ${\cal L}'_p = {\cal L}_p {\cal W}^*_p$ and denote the elements of 
the row ${\cal L}'_p = (L'_{p,0},\dots,L'_{p,\#{\cal E}_p-1})$.

\begin{theo}
	Let a dilation matrix $M$ and a center
	$c$ be  appropriate for an abelian symmetry group ${\cal H}$, $n\in\n$.
	Suppose $m_0$ and $\w m_0$ are as in Theorem~\ref{theoWaveSym} and 
	condition~(\ref{fSpecAssum}) is valid for all 
	$p=0,\dots,\#\Lambda-1.$
	Suppose $m_{(p,i)},$ $\w m_{(p,i)}$ are wavelet masks constructed using
	Theorem~\ref{theoWaveSym}.
	Assume that trigonometric polynomials $L_{p,i}$ are as in Theorem~\ref{theoLiftSym}.	
	New masks $m_{(p,i)}^n$, $\w m_{(p,i)}^n$ defined by the lifting
	scheme transformation with trigonometric polynomials $L'_{p,i}$ have the same symmetry properties as masks $m_{(p,i)}$, $\w m_{(p,i)}$.
	New wavelet masks $m_{(p,i)}^n$ have vanishing moments at least of order 1.
	If new refinable functions $\phi,$ $\w\phi$
	corresponding to new refinable masks $m^n_0$, $\w m^n_0$
	are in $L_2(\rd)$, then the resulting wavelet system
	is a dual wavelet frame.
\end{theo}

{\bf Proof.} By direct computations we can state that
$$
L'_{p,i}((M^{-1}E^{(k)} M)^*\xi)=
{\left[{\cal W}_p\right]_{k,i}} L'_{p,i}(\xi), \quad \forall E^{(k)} \in {\cal E}_p
$$
and
$$
L'_{p,i}((M^{-1}F M)^*\xi)=
L'_{p,i}(\xi)e^{2\pi i (r_{p,0}^F,\xi)}\quad \forall F\in {\cal H}_{p,0},
$$
$i = 0,\dots,\#{\cal E}_p-1,$ $p=0,\dots,\#\Lambda-1$.
Therefore, we obtain
$$
m_{(p,i)}^n(F^*\xi) = e^{2\pi i (F s_{p,0} - s_{p,0},\xi)}
m_{(p,i)}^n(\xi)
$$
and 
$$
m_{(p,i)}^n(E^{(k)*}\xi) =[{\cal W}_p]_{k,i} e^{2\pi i (E^{(k)} c - c,\xi)}
m_{(p,i)}^n(\xi),
$$
$i = 0,\dots,\#{\cal E}_p-1,$ $p=0,\dots,\#\Lambda-1$.
Thus, wavelet mask $m_{(p,i)}^n$ has the same symmetry properties as $m_{(p,i)}$.
Also, by direct computations it can be checked that $\w m_0^n$ remains ${\cal H}$-symmetric with respect to the center $c$. $\Diamond$

\subsection{Symmetrization for frames}

Symmetrization step also can be done in Theorem~\ref{theoFrameNoSym}.
Again, we assume that ${\cal H}$ is an abelian symmetry group and~(\ref{fSpecAssum}) is valid for all $p=0,\dots,\#\Lambda-1$.
For $r=m$, instead of matrix extension~(\ref{eq:MEPframe1}) we consider
\be
{\cal Q}:=\left(%
\begin{array}{cc}
	P & 0 \\
	{\cal W}^* - {\cal W}^*\w P^* P & {\cal W}^* \w P^*\\
\end{array}%
\right),\ \ \
\widetilde{\cal Q}:=\left(%
\begin{array}{cc} 
	\w P & 0 \\
	{\cal W}^* - {\cal W}^* P^* \w P & {\cal W}^* P^*\\
\end{array}%
\right),
\label{eq:MEPframe3}
\ee
where  matrix ${\cal W}$ is a symmetrization matrix.
Note that ${\cal Q}^* \w{\cal Q} = I_r$, since 
${\cal Q} = \texttt{diag}(1,{\cal W}^*)\, {\cal N}$ and 
$\w {\cal Q} = \texttt{diag}(1,{\cal W}^*)\, \w {\cal N}$.

For $r=m+1$, instead of matrix extension~(\ref{eq:MEPframe2}) we consider
\be
{\cal N}:=\left(%
\begin{array}{ccc}
	P & \mu_{0m} & 0 \\
	{\cal W}^* - {\cal W}^* \w P^* P & \mu_{0m} {\cal W}^* P^* & {\cal W}^* \w P^*\\
	- \overline{\w\mu_{0m}} P &  1 - \overline{\w\mu_{0m}} \mu_{0m} & \overline{\w\mu_{0m}}\\
\end{array}%
\right),\ \ \
\widetilde{\cal N}:=\left(%
\begin{array}{ccc} 
	\w P & \w\mu_{0m} & 0 \\
	{\cal W}^* - {\cal W}^* P^* \w P & \w \mu_{0m} {\cal W}^* P^* & {\cal W}^* P^*\\
	- \overline{\mu_{0m}} \w P &  1 - \overline{\mu_{0m}} \w\mu_{0m} & \overline{\mu_{0m}}
\end{array}%
\right), 
\label{eq:MEPframe4}
\ee
Note that ${\cal Q}^* \w {\cal Q} = I_r$, since ${\cal Q}=\texttt{diag}(1,{\cal W}^*,1)\, {\cal N}$, 
$\widetilde{\cal Q}=\texttt{diag}(1,{\cal W}^*,1)\, \w{\cal N}$.

\begin{theo} 
	Let a dilation matrix $M$ and a center
	$c$ be  appropriate for an abelian symmetry group ${\cal H}$, $n\in\n$.
	Suppose masks $m_\nu, \w m_\nu$, $\nu=0\ddd r$ are constructed by
	Algorithm 1 where the matrix extension is realized by~(\ref{eq:MEPframe3}), if $r=m$, or by~(\ref{eq:MEPframe4}), if $r=m+1$.
	Then wavelet masks  $m_{(p,i)}$ and $\w m_{(p,i)}$
	have the ${\cal H}$-symmetry properties, 
	i.e. for 
	$K\in {\cal H}$, with 
	$K=E^{(k)}F$, $E^{(k)}\in {\cal E}_{p}$,
	$F\in {\cal H}_{p,0},$ 
	$$
	m_{(p,i)}(K^*\xi) = 
	[W_p]_{k,i}  e^{2\pi i (Kc - c + M r_{p,0}^F,\xi)} m_{(p,i)}(\xi)
	$$
	and analogously for $\w m_{(p,i)}$,
	$i=0,\dots,\#{\cal E}_{p}-1,$
	$p=0,\dots,\#\Lambda-1$.
	If $r=m+1$, then the last wavelet masks $m_r$ and $\w m_r$ defined by 
	the last rows of matrices ${\cal N}$ and $\w {\cal N}$ are 
	${\cal H}$-symmetric with respect to the center $c$.
	All wavelet masks  have vanishing moments of order  $n$.  If refinable functions $\phi,$ $\w\phi$
	corresponding to refinable masks $m_0$, $\w m_0$
	are in $L_2(\rd)$, then the resulting wavelet system
	is a dual wavelet frame.
	\label{theoFrameSym}
\end{theo}

The proof can be done analogously to the proof of Theorem~\ref{theoWaveSym}.

\section{Examples}

In this section we
give several examples which illustrate the results of the paper.	
For a compactly supported distribution $f\in {\cal S}',$ its
\emph{Sobolev smoothness exponent} is defined to be
$$\nu_2(f)=\sup \left\{ \nu\in{\mathbb R}\cup\{-\infty\}\cup\{+\infty\}: \int_{{\mathbb R}^d} |\widehat f (\xi)|^2 (1+ \|\xi\|_2^2)^{\nu} d \xi < \infty\right\}.$$
Below, the Sobolev smoothness exponent is calculated by Theorem 7.1 in~\cite{HanSymSmooth}.

1. Let ${\cal H}$ be a hexagonal symmetry group on $\z^2$, namely,
\begin{small}
	
	$${\cal H}=\left\{\pm I_2,\pm\left(
	\begin{array}{cc}
	0 & 1 \\
	1 & 0 \\
	\end{array}
	\right),\pm\left(
	\begin{array}{cc}
	1 & 0 \\
	1 & -1 \\
	\end{array}
	\right),\pm\left(
	\begin{array}{cc}
	1 & -1 \\
	1 & 0 \\
	\end{array}
	\right),\pm\left(
	\begin{array}{cc}
	0 & 1 \\
	-1 & 1 \\
	\end{array}
	\right),\pm\left(
	\begin{array}{cc}
	-1 & 1 \\
	0 & 1 \\
	\end{array}
	\right)\right\}$$
\end{small}
$c={\bf0},$ ${M}= \left(\begin{array}{cc} 2&0 \\
0& 2 
\end{array}\right).
$
The set of digits is $D(M)=\{s_0=(0,0), s_1=(0,1), s_2=(-1,0),s_3=(-1,-1)\},$ $m=4.$ 
Let us construct an interpolatory refinable mask  that is ${\cal H}$-symmetric with respect to the origin
and obeys the sum rule of order $n=3.$
The digits are renumbered as $s_{0,0}=(0,0),$ $s_{1,0}=(0,1),s_{1,1}=(-1,0), s_{1,2} = (-1,-1).$
And ${\cal H}_{0,0}={\cal H},$ ${\cal E}_{0}=\{I_2\},$
$${\cal H}_{1,0}=\left\{\pm I_2, \pm \left(
\begin{array}{cc}
1 & -1 \\
1 & 0 \\
\end{array}
\right)\right\},$$
${\cal E}_{1}=\left\{ E^{(0)}, E^{(1)},E^{(2)}\right\} = \left\{I_2,\left(
\begin{array}{cc}
1 & -1 \\
1 & 0 \\
\end{array}
\right), \left(
\begin{array}{cc}
1 & -1 \\
0 & -1 \\
\end{array}
\right)\right\}.$
According to Theorem 10 in~\cite{KrSymInterp},
mask $m_0$
can be constructed as follows

$$m_0: \left(
\begin{array}{ccccccc}
0 & 0 & 0 & -\frac{1}{64} & 0 & 0 & -\frac{1}{64} \\
0 & 0 & 0 & 0 & 0 & 0 & 0 \\
0 & 0 & 0 & \frac{9}{64} & \frac{9}{64} & 0 & 0 \\
-\frac{1}{64} & 0 & \frac{9}{64} & \frac{1}{4} & \frac{9}{64} & 0 & -\frac{1}{64} \\
0 & 0 & \frac{9}{64} & \frac{9}{64} & 0 & 0 & 0 \\
0 & 0 & 0 & 0 & 0 & 0 & 0 \\
-\frac{1}{64} & 0 & 0 & -\frac{1}{64} & 0 & 0 & 0 \\
\end{array}
\right)$$
with coefficient support in $[-3,3]^2\bigcap{\mathbb Z}^2.$ Note that $\mu_{00}=\frac 12.$
The corresponding refinable function $\phi$ is in $L_2({\mathbb R}^2)$ since $\nu_2(\phi)\ge 1.76585.$
Dual mask can be taken as $\w m_0 \equiv 1.$ Then by Theorem~\ref{theoWaveNoSym} we get wavelet masks
$$
m_{1,0}: \left(
\begin{array}{ccccccc}
0 & 0 & 0 & \frac{1}{512} & 0 & 0 & \frac{1}{512} \\
0 & 0 & 0 & 0 & 0 & 0 & 0 \\
0 & 0 & 0 & -\frac{1}{64} & -\frac{9}{512} & 0 & \frac{1}{512} \\
\frac{1}{512} & 0 & -\frac{9}{512} & -\frac{1}{32} & -\frac{9}{512} & 0 & \frac{1}{512} \\
0 & 0 & -\frac{9}{512} & \frac{55}{256} & -\frac{9}{512} & 0 & 0 \\
\frac{1}{512} & 0 & -\frac{9}{512} & -\frac{1}{32} & -\frac{9}{512} & 0 & \frac{1}{512} \\
\frac{1}{512} & 0 & -\frac{9}{512} & -\frac{1}{64} & 0 & 0 & 0 \\
0 & 0 & 0 & 0 & 0 & 0 & 0 \\
\frac{1}{512} & 0 & 0 & \frac{1}{512} & 0 & 0 & 0 \\
\end{array}
\right),$$
with coefficient support in $[-3,3]\times[-4,4]\bigcap{\mathbb Z}^2,$
$ m_{1,1}(\xi) = m_{1,0} (E^{(1)*}\xi), m_{1,2}(\xi) =  m_{1,0}  (E^{(2)*}\xi).
$
The dual wavelet masks are $\w m_{1,0} (\xi) = \frac 14 e^{2\pi i \xi_2},$ $\w m_{1,1} (\xi)= \frac 14 e^{-2\pi i \xi_1},$ $\w m_{1,2}(\xi) =\frac 14 e^{-2\pi i (\xi_1+\xi_2)}$.
Since the symmetry group is not abelian, the symmetrization step cannot be done here. Thus, we get mutually symmetric frame-like wavelet system providing approximation order $3$. 
Now, using the lifting scheme, we improve this system to dual wavelet frames keeping the symmetry and approximation properties.
Let $L_1(\xi) = -\frac 18 - \frac 18 e^{2\pi i \xi_2}$, $L_2 (\xi) = -\frac 18 - \frac 18 e^{-2\pi i \xi_1},$ $L_3(\xi) = -\frac 18 - \frac 18 e^{-2\pi i (\xi_1+\xi_2)}$. Then the new dual wavelet masks are $$\w m_{1,0} (\xi) = e^{2\pi i \xi_2} \left(1 - \frac{9}{16} e^{2\pi i \xi_2} - \frac{9}{16} e^{-2\pi i \xi_2} + \frac{1}{16} e^{2\pi i 3 \xi_2} +  \frac{1}{16} e^{-2\pi i 3 \xi_2}\right),$$ $ m_{1,1}(\xi) = m_{1,0} (E^{(1)*}\xi), m_{1,2}(\xi) =  m_{1,0}  (E^{(2)*}\xi).
$
The new  dual refinable mask is 
$$
\w m_0: \left(
\begin{array}{ccccccccc}
0 & 0 & 0 & 0 & \frac{1}{128} & 0 & 0 & 0 & \frac{1}{128} \\
0 & 0 & 0 & 0 & 0 & 0 & 0 & 0 & 0 \\
0 & 0 & 0 & 0 & -\frac{1}{16} & 0 & -\frac{1}{16} & 0 & 0 \\
0 & 0 & 0 & 0 & \frac{1}{8} & \frac{1}{8} & 0 & 0 & 0 \\
\frac{1}{128} & 0 & -\frac{1}{16} & \frac{1}{8} & \frac{37}{64} & \frac{1}{8} & -\frac{1}{16} & 0 & \frac{1}{128} \\
0 & 0 & 0 & \frac{1}{8} & \frac{1}{8} & 0 & 0 & 0 & 0 \\
0 & 0 & -\frac{1}{16} & 0 & -\frac{1}{16} & 0 & 0 & 0 & 0 \\
0 & 0 & 0 & 0 & 0 & 0 & 0 & 0 & 0 \\
\frac{1}{128} & 0 & 0 & 0 & \frac{1}{128} & 0 & 0 & 0 & 0 \\
\end{array}
\right)
$$
with coefficient support in $[-4,4]^2\bigcap{\mathbb Z}^2$. 
Since $\w\phi$ is in $L_2(\rd)$ ($\nu_2(\phi)\ge 0.1566$), we get mutually symmetric dual wavelet frame in $L_2(\rd)$ providing approximation order $3$.

2. 1. Let ${\cal H}$ be a point symmetry group on $\z^2$, namely,	
${\cal H}=\left\{\pm I_2\right\},$
$c=(1/2,0),$ ${M}= \left(\begin{array}{cc} 1&-2\cr
2& -1 
\end{array}\right).
$
The set of digits is $D(M)=\{s_0=(0,0), s_1=(-1,0), s_2=(1,0)\},$ $m=3.$ 
Let us construct a refinable mask  that is ${\cal H}$-symmetric with respect to $c$
and obeys the sum rule of order $n=2.$
The digits are renumbered as $s_{0,0}=(-1,0),$ $s_{1,0}=(0,0),s_{1,1}=(1,0).$
And ${\cal H}_{0,0}={\cal H},$ ${\cal E}_{0}=\{I_2\},$
${\cal H}_{1,0}=\left\{ I_2\right\},$
${\cal E}_{1}= {\cal H}$.
According to Theorem 10 in~\cite{KrSymInterp},
mask $m_0$
can be constructed as follows

$$m_0: \left(
\begin{array}{cccc}
0 & -\frac{1}{12} & 0 & \frac{1}{18} \\
0 & 0 & \frac{1}{6} & \frac{1}{9} \\
\frac{1}{12} & \frac{1}{6} & \frac{1}{6} & \frac{1}{12} \\
\frac{1}{9} & \frac{1}{6} & 0 & 0 \\
\frac{1}{18} & 0 & -\frac{1}{12} & 0 \\
\end{array}
\right)$$
with coefficient support in $[-1,2]\times [-2,2]\bigcap{\mathbb Z}^2.$ 
The corresponding refinable function $\phi$ is in $L_2({\mathbb R}^2)$ since $\nu_2(\phi)\ge 0.776.$
The utility dual mask can be taken as $\w m'_0(\xi) = 1/3 + 1/3 e^{2\pi i \xi_1} + 1/6 e^{-2\pi i \xi_2} + 1/6 e^{2\pi i (\xi_1+\xi_2)}.$ This mask obeys the sum rule of order $1$. Then $\sigma = \sum\limits_{l=0}^{m-1}\overline{\mu_{0l}}\,\widetilde\mu'_{0l}$ satisfied $D^{\beta}(1-\sigma)(\nul) = 0$ for $|\beta|<n$. Then by Algorithm 1 we get dual mask
$$
\w m_0:  \left(
\begin{array}{cccccc}
0 & 0 & 0 & \frac{1}{144} & 0 & 0 \\
0 & 0 & \frac{1}{72} & \frac{1}{72} & -\frac{5}{216} & 0 \\
0 & 0 & \frac{1}{144} & -\frac{5}{108} & -\frac{5}{108} & -\frac{11}{432} \\
0 & 0 & 0 & \frac{49}{216} & -\frac{11}{216} & -\frac{11}{216} \\
0 & -\frac{11}{432} & \frac{1}{2} & \frac{1}{2} & -\frac{11}{432} & 0 \\
-\frac{11}{216} & -\frac{11}{216} & \frac{49}{216} & 0 & 0 & 0 \\
-\frac{11}{432} & -\frac{5}{108} & -\frac{5}{108} & \frac{1}{144} & 0 & 0 \\
0 & -\frac{5}{216} & \frac{1}{72} & \frac{1}{72} & 0 & 0 \\
0 & 0 & \frac{1}{144} & 0 & 0 & 0 \\
\end{array}
\right)
$$
with coefficient support in $[-2,3]\times [-3,3]\bigcap{\mathbb Z}^2.$ 
Note that $\w\phi \in L_2 (\rd)$ since $\nu_2(\w\phi)\ge 0.503.$ Also,  $\w m_0$ obeys the sum rule of order $1$. Wavelet function and dual wavelet functions can be  constructed by the technique in Algorithm 1. 
Then wavelet functions have vanishing moments  of order $1$, dual wavelet functions have vanishing moments of order $2$.




\end{document}